\documentclass[amstex,12pt, amssymb]{article}

\usepackage{mathtext}
\usepackage[cp1251]{inputenc}
\usepackage[T2A]{fontenc}
\usepackage[dvips]{graphicx}
\usepackage{amsmath}
\usepackage{amssymb}
\usepackage{amsxtra}
\usepackage{latexsym}
\usepackage{ifthen}

\newcounter{lemma}[section]

\newcounter{corollary}[section]

\newcounter{remark}[section]

\newcounter{theorem}[section]

\newcounter{proposition}[section]

\numberwithin{equation}{section}

\def\Xint#1{\mathchoice
   {\XXint\displaystyle\textstyle{#1}}%
   {\XXint\textstyle\scriptstyle{#1}}%
   {\XXint\scriptstyle\scriptscriptstyle{#1}}%
   {\XXint\scriptscriptstyle\scriptscriptstyle{#1}}%
   \!\int}
\def\XXint#1#2#3{{\setbox0=\hbox{$#1{#2#3}{\int}$}
     \vcenter{\hbox{$#2#3$}}\kern-.5\wd0}}
\def\dashint{\Xint-}

\def\Rn{{{\Bbb R}^n}}
\def\Rk{{{\Bbb R}^k}}

\def\lRn{{\overline {{\Bbb R}^n}}}

\def\cc{\setcounter{equation}{0}
\setcounter{figure}{0}\setcounter{table}{0}}

\begin{document}

\markboth{\centerline{DENIS KOVTONYUK AND VLADIMIR RYAZANOV}}
{\centerline{THE THEORY OF PRIME ENDS AND SPATIAL MAPPINGS IV}}

\author{{DENIS KOVTONYUK AND VLADIMIR RYAZANOV}}

\title{{\bf THE THEORY OF PRIME ENDS \\ AND SPATIAL MAPPINGS IV}}

\maketitle

\large \begin{abstract} It is given a canonical representation of
prime ends in regular spatial domains and, on this basis, it is
studied the boundary behavior of the so-called lower
$Q$-homeomorphisms that are the natural generalization of the
quasiconformal mappings. In particular, it is found a series of
effective conditions on the function $Q(x)$ for a homeomorphic
extension of the given mappings to the boundary by prime ends in
domains with regular boundaries. The developed theory is applied, in
particular, to mappings of the classes of Sobolev and
Orlicz--Sobolev and also to finitely bi--Lipschitz mappings that a
far--reaching extension of the well--known classes of isometric and
quasiisometric mappings.
\end{abstract}

\bigskip
{\bf 2010 Mathematics Subject Classification: Primary 30C65,
30C85,30D40, 31A15, 31A20, 31A25, 31B25. Se\-con\-da\-ry 37E30.}

\large \cc
\section{Introduction}

The problem of the boundary behavior is one of the central topics of
the theory of quasiconformal mappings and their generalizations.
During the last years they intensively studied various classes of
mappings with finite distortion in a natural way generalizing
conformal, quasiconformal and quasiregular mappings, see many
references in the monographs \cite{GRSY} and \cite{MRSY}. In this
case, as it was earlier, the main geometric approach in the modern
mapping theory is the method of moduli, see, e.g., the monographs
\cite{GRSY}, \cite{MRSY}, \cite{Oht}, \cite{Ri}, \cite{Vs},
\cite{Va}  and \cite{Vu}.

From the point of view of the theory of conformal mappings, it was
un\-sa\-tis\-fac\-to\-ry to consider the individual points of the
boundary of a simply connected domain as the primitive constituents
of the boundary. Indeed, if correspondingly to the Riemann theorem
such a domain is mapped conformally onto the unit disk, then the
points of the unit circumference correspond to the so--called prime
ends of the domain.

The term "prime end"\ originated from Caratheodory \cite{Car$_2$}
who initiated the systematic study of the structure of the boundary
of a simply connected domain. His approach was topological and dealt
with concepts subdomains, crosscuts etc. that are defined with
reference to the given domain. The problem arisen under his approach
to show that prime ends are preserved under conformal mappings was
just solved by one of Caratheodory's fundamental theorems.

Lindel\"of \cite{Lind} circumvented this difficulty by defining
prime ends of a domain with reference to the conformal map of the
unit disk onto the domain; namely in terms of the set of
indetermination or cluster set. However, his method does not obviate
an explicit analysis of the topological situation in the domain
itself.

Two other schemes for the definition of prime ends deserve brief
mention. Mazurkiewicz \cite{Ma} introduced a metric
$\rho_{\pi}(z_1,z_2)$ that is equivalent to the euc\-li\-dean metric
in a domain in the sense that $\rho_{\pi}(z_j,z_0)\to0$ if and only
if $|z_j-z_0|\to0$ for any sequence $\{z_j\}$ of points of the
domain. The boundary of the domain with respect to $\rho_{\pi}$,
i.e. the complement of the domain with respect to its
$\rho_{\pi}-$completion, is a space that can be identified with the
set of prime ends of Caratheodory.

Finally, Ursell and Young \cite{UY} to introduce the prime ends of a
domain have used the notion of an equivalence class of paths that
converge to the boundary of the domain. For the history of the
question, see also \cite{ABBS}, \cite{CL} and \cite{Na} and further
references therein.

Later on, we often use the notations $\mathrm I$, $\bar{\mathrm I}$,
$\mathbb R$, $\overline{\mathbb R}$, ${\Bbb R}^+$, $\overline{{\Bbb
R}^+}$ and $\overline{{\Bbb R}^n}$ for $[1,\infty)$, $[1,\infty]$,
$(-\infty,\infty)$, $[-\infty ,\infty ]$, $[0,\infty)$, $[0,\infty
]$ and ${\Bbb R}^n\cup\{\infty\}$, correspondingly, and $D$ is a
domain in ${\Bbb R}^n$.

In what follows, we use in $\overline{{{\Bbb R}}^n}$ the {\bf
spherical (chordal) metric} $h(x,y)=|\pi(x)-\pi(y)|$ where $\pi$ is
the stereographic projection of $\overline{{{\Bbb R}}^n}$ onto the
sphere $S^n(\frac{1}{2}e_{n+1},\frac{1}{2})$ in ${{\Bbb R}}^{n+1},$
i.e.
$$h(x,y)=\frac{|x-y|}{\sqrt{1+{|x|}^2} \sqrt{1+{|y|}^2}},\ \,\, x\ne
\infty\ne y, \ \ \ \ \ \ h(x,\infty)=\frac{1}{\sqrt{1+{|x|}^2}}\ .
$$
The quantity
$$h(E)=\sup_{x,y \in E} h(x,y)$$
is said to be {\bf spherical (chordal) diameter} of a set $E \subset
\overline{{\Bbb R}^n}$.

Let $\omega$ be an open set in $\Rk$, $k=1,\ldots,n-1$. A
(continuous) mapping $\sigma:\omega\to\lRn$ is called a {\bf
$k-$dimensional surface} in $\lRn$. An $(n-1)-$dimensional surface
$\sigma$ in $\lRn$ is called also a {\bf surface}. A surface
$\sigma:\omega\to D$ is called a {\bf Jordan surface in} $D$ if
$\sigma(z_1)\neq\sigma(z_2)$ whenever $z_1\neq z_2$. Later on, we
sometimes use $\sigma$ to denote the whole image
$\sigma(\omega)\subseteq\lRn$ under the mapping $\sigma$, and
$\overline{\sigma}$ instead of $\overline{\sigma(\omega)}$ in $\lRn$
and $\partial\sigma$ instead of
$\overline{\sigma(\omega)}\setminus\sigma(\omega)$. A Jordan surface
$\sigma$ in $D$ is called a {\bf cut} of $D$ if $\sigma$ splits $D$,
i.e. $D\setminus \sigma$ has more than one component,
$\partial\sigma\cap D=\varnothing$ and $\partial\sigma\cap\partial
D\neq\varnothing$.
\medskip

A sequence $\sigma_1,\ldots, \sigma_m,\ldots$ of cross--cuts of $D$
is called a {\bf chain} if:
\medskip

(i) $\overline{\sigma_i}\cap\overline{\sigma_j}=\varnothing$ for
every $i\neq j$, $i,j= 1,2,\ldots$;
\medskip

(ii) $\sigma_{m-1}$ and $\sigma_{m+1}$ are contained in different
components of $D\setminus \sigma_m$ for every $m>1$;
\medskip

(iii) $\cap\,d_m=\varnothing$ where $d_m$ is a component of
$D\setminus \sigma_m$ containing $\sigma_{m+1}$.

\bigskip

Finally, we will call a chain of cross--cuts $\{\sigma_m\}$ {\bf
regular} if

\medskip

(iv)  $h(\sigma_{m})\to0$ as $m\to\infty$.

\bigskip

Correspondingly to the definition, a chain of cross--cuts
$\{\sigma_m\}$ is determined by a chain of domains $d_m\subset D$
such that $\partial\,d_m\cap D\subseteq\sigma_m$ and $d_1\supset
d_2\supset\ldots\supset d_m\supset\ldots$. Two chains of cross--cuts
$\{\sigma_m\}$ and $\{\sigma_k'\}$ are called {\bf equivalent} if,
for every $m=1,2,\ldots$, the domain $d_m$ contains all domains
$d_k'$ except its finite collection and, for every $k=1,2,\ldots$,
the domain $d_k'$ contains all domains $d_m$ except its finite
collection, too. An {\bf end} $K$ of the domain $D$ is an
equivalence class of chains of cross--cuts of $D$.

Let $K$ be an end of a domain $D$ in $\lRn$ and $\{\sigma_m\}$ and
$\{\sigma_m'\}$ be two chains in $K$ and $d_m$ and $d_m'$  be
domains corresponding to $\sigma_m$ and $\sigma_m'$, respectively.
Then
$$\bigcap\limits_{m=1}\limits^{\infty}\overline{d_m}\ \subseteq\
\bigcap\limits_{m=1}\limits^{\infty}\overline{d_m'}\ \subset\
\bigcap\limits_{m=1}\limits^{\infty}\overline{d_m}$$ and, thus,
$$\bigcap\limits_{m=1}\limits^{\infty}\overline{d_m}\ =\ \bigcap\limits_{m=1}\limits^{\infty}\overline{d_m'}\ ,$$
i.e. the set
$$I(K)\ =\ \bigcap\limits_{m=1}\limits^{\infty}\overline{d_m}$$ depends only on $K$ but not on
a choice of its chain of cross--cuts $\{\sigma_m\}$. The set $I(K)$
is called the {\bf impression of the end} $K$. It is well--known
that $I(K)$ is a continuum, i.e. a connected compact set, see, e.g.,
I(9.12) in \cite{Wh}. Moreover, in view of the conditions (ii) and
(iii), we obtain that
$$I(K)\ =\ \bigcap\limits_{m=1}\limits^{\infty}(\partial d_m\cap\partial D)\ =\
\partial D\ \cap\ \bigcap\limits_{m=1}\limits^{\infty}\partial d_m\,.$$
Thus, we come to the following conclusion.

\bigskip

\begin{proposition}\label{thabc1} {\it For every end $K$  of a domain $D$
in $\overline{{\Bbb R}^n}$,}
\begin{equation}\label{eq1}I(K)\subseteq\partial D\,.\end{equation}
\end{proposition}

Following \cite{Na}, we say that $K$ is a {\bf prime end} if $K$
contains a chain of cross--cuts $\{\sigma_m\}$ such that
\begin{equation}\label{eqSIMPLE}
\lim\limits_{m\to\infty}\ M(\Delta(C,\sigma_m;D))\ =\ 0
\end{equation} for a continuum $C$ in $D$ where $\Delta(C,\sigma_m;D)$
is the collection of all paths connecting the sets $C$ and
$\sigma_m$ in $D$ and $M$ denotes its modulus, see the next section.

If an end $K$ contains at least one regular chain, then $K$ will be
said to be {\bf regular}. As it will easy follow from Lemma
\ref{thabc2}, every regular end is a prime end.

\cc
\section{On lower and ring $Q$-homeomorphisms}

The class of lower $Q$-homeomorphisms was introduced in the paper
\cite{KR$_1$}, see also the monograph \cite{MRSY}, and was motivated
by the ring definition of quasiconformal mappings of Gehring, see
\cite{Ge$_1$}. The theory of lower $Q$-homeomorphisms has found
interesting applications to the theory of the Beltrami equations in
the plane and to the theory of mappings of the classes of Sobolev
and Orlich-Sobolev in the space, see, e.g., \cite{KPR},
 \cite{KPRS}, \cite{KRSS1}, \cite{KRSS}, \cite{KSS}, \cite{MRSY} and
\cite{RSSY}.

\medskip

Let $\omega$ be an open set in $\overline{{\Bbb R}^k}$,
$k=1,\ldots,n-1$. Recall that a (continuous) mapping
$S:\omega\to{\Bbb R}^n$ is called a $k$-dimensional surface $S$ in
${\Bbb R}^n$. The number of preimages
\begin{equation}\label{eq8.2.3} N(S,y)={\rm card}\,S^{-1}(y)={\rm
card}\,\{x\in\omega:S(x)=y\},\ y\in{\Bbb R}^n\end{equation} is said
to be a {\bf multiplicity function} of the surface $S$. It is known
that the multiplicity function is lower semicontinuous, i.e.,
$$N(S,y)\geqslant\liminf_{m\to\infty}\:N(S,y_m)$$ for every sequence $y_m\in{\Bbb R}^n$, $m=1,2,\ldots\,$, such
that $y_m\to y\in{\Bbb R}^n$ as $m\to\infty$, see, e.g., \cite{RR},
p. 160. Thus, the function $N(S,y)$ is Borel measurable and hence
measurable with respect to every Hausdorff measure $H^k$, see, e.g.,
\cite{Sa}, p. 52.

\medskip

Recall that a $k$-dimensional Hausdorff area in ${\Bbb R}^n$ (or
simply {\bf area}) associated with a surface $S:\omega\to{\Bbb R}^n$
is given by \begin{equation}\label{eq8.2.4}
{\mathcal{A}}_S(B)={\mathcal{A}}^{k}_S(B):=\int\limits_B
N(S,y)\,dH^{k}y
\end{equation} for every Borel set $B\subseteq{\Bbb R}^n$ and,
more generally, for an arbitrary set that is measurable with respect
to $H^k$ in ${\Bbb R}^n$, cf. 3.2.1 in \cite{Fe} and 9.2 in
\cite{MRSY}.

\medskip

If $\varrho:{\Bbb R}^n\to\overline{{\Bbb R}^+}$ is a Borel function,
then its {\bf integral over} $S$ is defined by the equality
\begin{equation}\label{eq8.2.5} \int\limits_S \varrho\,d{\mathcal{A}}:=
\int\limits_{{\Bbb R}^n}\varrho(y)\,N(S,y)\,dH^ky\,.\end{equation}
Given a family $\Gamma$ of $k$-dimensional surfaces $S$, a Borel
function $\varrho:{\Bbb R}^n\to[0,\infty]$ is called {\bf
admissible} for $\Gamma$, abbr. $\varrho\in{\rm adm}\,\Gamma$, if
\begin{equation}\label{eq8.2.6}\int\limits_S\varrho^k\,d{\mathcal{A}}\geqslant1\end{equation}
for every $S\in\Gamma$. The {\bf modulus} of $\Gamma$ is the
quantity
\begin{equation}\label{eq8.2.7} M(\Gamma)=\inf_{\varrho\in{\rm adm}\,\Gamma}
\int\limits_{{\Bbb R}^n}\varrho^n(x)\,dm(x)\,.\end{equation} We also
say that a Lebesgue measurable function $\varrho:{\Bbb
R}^n\to[0,\infty]$ is {\bf extensively admissible} for a family
$\Gamma$ of $k$-dimensional surfaces $S$ in ${\Bbb R}^n$, abbr.
$\varrho\in{\rm ext\, adm}\,\Gamma$,  if a subfamily of all surfaces
$S$ in $\Gamma$, for which (\ref{eq8.2.6}) fails, has the modulus
zero.

Given domains $D$ and $D'$ in $\overline{{\Bbb R}^n}={\Bbb
R}^n\cup\{\infty\}$, $n\geqslant2$,
$x_0\in\overline{D}\setminus\{\infty\}$, and a measurable function
$Q:{{\Bbb R}^n}\to(0,\infty)$, we say that a homeomorphism $f:D\to
D'$ is a {\bf lower $Q$-homeomorphism at the point} $x_0$ if
\begin{equation}\label{eqOS1.10} M(f\Sigma_{\varepsilon})\geqslant
\inf\limits_{\varrho\in{\rm
ext\,adm}\,\Sigma_{\varepsilon}}\int\limits_{D\cap
R_{\varepsilon}}\frac{\varrho^n(x)}{Q(x)}\,dm(x)\end{equation} for
every ring $R_{\varepsilon}=\{x\in{\Bbb
R}^n:\varepsilon<|x-x_0|<\varepsilon_0\}\,,\quad\varepsilon\in(0,\varepsilon_0)\,,\
\varepsilon_0\in(0,d_0)$, where $d_0=\sup\limits_{x\in D}|x-x_0|$,
and $\Sigma_{\varepsilon}$ denotes the family of all intersections
of the spheres $S(x_0,r)=\{x\in{\Bbb R}^n:|x-x_0|=r\}\,,
r\in(\varepsilon,\varepsilon_0)\,,$ with $D$. This notion can be
extended to the case $x_0=\infty\in\overline{D}$ by applying the
inversion $T$ with respect to the unit sphere in $\overline{{\Bbb
R}^n}$, $T(x)=x/|x|^2$, $T(\infty)=0$, $T(0)=\infty$. Namely, a
homeomorphism $f:D\to D'$ is said to be a {\bf lower
$Q$-homeomorphism at} $\infty\in\overline{D}$ if $F=f\circ T$ is a
lower $Q_*$-homeomorphism with $Q_*=Q\circ T$ at $0$.

\medskip

We also say that a homeomorphism $f:D\to{\overline{{\Bbb R}^n}}$ is
a {\bf lower $Q$-homeomor\-phism in} $D$ if $f$ is a lower
$Q$-homeomorphism at every point $x_0\in\overline{D}$.

\medskip

Recall the criterion for homeomorphisms in ${\Bbb R}^n$ to be lower
$Q$-homeomorphisms, see Theorem 2.1 in \cite{KR$_1$} or Theorem 9.2
in \cite{MRSY}.

\medskip

\begin{proposition}\label{prOS2.2}
{\it Let $D$ and $D'$ be domains in $\overline{{\Bbb R}^n}$,
$n\geqslant2$, let $x_0\in\overline{D}\setminus\{\infty\}$, and
$Q:D\to(0,\infty)$ be a measurable function. A homeomorphism $f:D\to
D'$ is a lower $Q$-homeomorphism at $x_0$ if and only if
\begin{equation}\label{eqOS2.1} M(f\Sigma_{\varepsilon})\geqslant\int\limits_{\varepsilon}^{\varepsilon_0}
\frac{dr}{||\,Q||_{n-1}(x_0,r)}\quad\quad\forall\
\varepsilon\in(0,\varepsilon_0)\,,\
\varepsilon_0\in(0,d(x_0))\,,\end{equation} where $d(x_0) =
\sup\limits_{x\in D}\,|\,x-x_0| $ and
\begin{equation}\label{eqOS2.3} ||Q||_{n-1}(x_0,r)=\left(\int\limits_{D(x_0,r)}Q^{n-1}(x)\,d{\mathcal A}\right)^\frac{1}{n-1}\end{equation} is the $L_{n-1}$-norm of $Q$ over $D(x_0,r)=\{x\in D:|x-x_0|=r\}=D\cap
S(x_0,r)$.}
\end{proposition}

\medskip

Further, as usual for sets $A$, $B$ and $C$ in $\overline{{\Bbb
R}^n}$, $\Delta(A,B,C)$ denotes the family of all paths joining $A$
and $B$ in $C$.

\medskip

Now, given domains $D$ in ${\Bbb R}^n$ and $D'$ in $\overline{{\Bbb
R}^n}$, $n\geqslant2$, and a measurable function $Q:{\Bbb
R}^n\to(0,\infty)$. Let $S_i:=S(x_0,r_i)$. We say that a
homeomorphism $f:D\to D'$ is a {\bf ring $Q$-homeomorphism at a
point} $x_0\in \overline{D}\setminus\{\infty\}$ if
\begin{equation}\label{eqOS1.8} M(f(\Delta(S_1,S_2,D)))\leqslant
\int\limits_{A\cap D}Q(x)\cdot\eta^{n}(|x-x_0|)\,dm(x)\end{equation}
for every ring $A=A(x_0,r_1,r_2)$, $0<r_1<r_2<d_0={\rm
dist}(x_0,\partial D)$, and for every measurable function
$\eta:(r_1,r_2)\to[0,\infty]$ such that
\begin{equation}\label{eqOS1.9}\int\limits_{r_1}^{r_2}\eta(r)\,dr\geq1.\end{equation}
The notion of a ring $Q-$homeomorphim can be extended to $\infty$ by
the standard way as in the case of a lower $Q-$homeomorphism above.

\medskip

The notion of a ring $Q-$homeomorphim was first introduced for inner
points of a domain in the work \cite{RSY$_1$} in the connection with
investigations of the Beltrami equations in the plane and then it
was extended to the space in the work \cite{RS}, see also the
monograph \cite{MRSY}. This notion was extended to boundary points
in the papers \cite{Lom} and \cite{RSY$_6$}--\cite{RSY$_2$}, see
also the monograph \cite{GRSY}. By Corollary 5 in \cite{KRSS} we
have the following fact.

\medskip

\begin{proposition}\label{corOS2.1} {\it In ${\Bbb R}^n$, $n\geqslant2$,
a lower $Q$-homeomorphism $f:D\to D'$ at a point $x_0\in
\overline{D}$ with $Q$ that is integrable in the degree $n-1$ in a
neighborhood of $x_0$ is a ring $Q_*$-homeomorphism at $x_0$ with
$Q_*=Q^{n-1}$.}
\end{proposition}

\medskip

\begin{remark}\label{rem7.4.3} By Remark 8 in \cite{KRSS} the conclusion
of Proposition \ref{corOS2.1} is valid if the function $Q$ is only
integrable in the degree $n-1$ on almost all spheres of small enough
radii centered at the point $x_0$.

Note also that, in the definitions of lower and ring
$Q-$homeomorphisms, it is sufficient to give the function $Q$ only
in the domain $D$ or to extend by zero outside of $D$.
\end{remark}

\cc
\section{On canonical representation of ends of spatial domains}

\begin{lemma}\label{thabc2} {\it Every regular end $K$ of a domain $D$ in $\lRn$
includes a chain of cross--cuts $\sigma_m$ lying on the spheres
$S_m$ centered at a point $x_0\in\partial D$ with hordal radii
$\rho_m\to0$ as $m\to\infty$. Every regular end $K$ of a bounded
domain $D$ in $\Rn$ includes a chain of cross--cuts $\sigma_m$ lying
on the spheres $S_m$ centered at a point $x_0\in\partial D$ with
euclidean radii $r_m\to0$ as $m\to\infty$.}
\end{lemma}

\medskip

{\bf Proof.} We restrict ourselves to the case of a domain $D$ in
$\lRn$ with the hordal metric. The second case is similar.

Let $\{\sigma_m\}$ be a chain of cross--cuts in the end $P$ and
$x_m$ a sequence of points in $\sigma_m$. Without loss of generality
we may assume that $x_m\to x_0\in\partial D$ as $m\to\infty$ because
$\lRn$ is a compact metric space. Then
$\rho^-_m:=h(x_0,\sigma_m)\to0$ because $h(\sigma_m)\to0$ as
$m\to\infty$. Furthermore, $$\rho^+_m\ :=\ H(x_0,\sigma_m)\ =\
\sup\limits_{x\in\sigma_m}h(x,x_0)\ =\
\sup\limits_{x\in\overline{\sigma_m}}h(x,x_0)$$ is the Hausdorff
distance between the compact sets $\{x_0\}$ and
$\overline{\sigma_m}$ in $\lRn$. By the condition (i) in the
definition of an end, we may assume without loss of generality that
$\rho^-_m>0$ and $\rho^+_{m+1}<\rho^-_m$ for all $m=1,2,\ldots$.

Set $$\delta_m\ =\ \Delta_m\setminus d_{m+1}$$ where
$\Delta_m=S_m\cap d_m$ and $$S_m\ =\ \{\ x\in\lRn\ :\
h(x_0,x)=\frac{1}{2}\left(\rho^-_m+\rho^+_{m+1}\right)\}\,.$$ It is
clear that $\Delta_m$ and $\delta_m$ are relatively closed in $d_m$.

Note that $d_{m+1}$ is contained in one of the components of the
open set $d_m\setminus\delta_m$. Indeed, assume that there is a pair
of points $x_1$ and $x_2\in d_{m+1}$ in different components
$\Omega_1$ and $\Omega_2$ of $d_m\setminus\delta_m$. Then $x_1$ and
$x_2$ can be joined by a continuous curve $\gamma:[0,1]\to d_{m+1}$.
However, $d_{m+1}$, and hence $\gamma$, does not intersect
$\delta_m$ by the construction and, consequently,
$[0,1]=\bigcup\limits_{k=1}^{\infty}\omega_k$ where
$\omega_k=\gamma^{-1}(\Omega_k)$, $\Omega_k$ is enumeration of
components $d_m\setminus\delta_m$. But $\omega_k$ are open in
$[0,1]$ because $\Omega_k$ are open and $\gamma$ is continuous. The
later contradicts to the connectivity of $[0,1]$ because
$\omega_1\neq\varnothing$ and $\omega_2\neq\varnothing$ and,
moreover, $\omega_i$ and $\omega_j$ are mutually disjoint whenever
$i\neq j$.

Let $d^*_m$ be a component of $d_m\setminus\delta_m$ containing
$d_{m+1}$. Then by the construction $d_{m+1}\subseteq d^*_m\subseteq
d_m$. It remains to show that $\partial d^*_m\setminus\partial
D\subseteq\delta_m$. First, it is clear that $\partial
d^*_m\setminus\partial D\subseteq\delta_m\cup\sigma_m$ because every
point in $d_m\setminus\delta_m$ belongs either to $d^*_m$ or to
other component of $d_m\setminus\delta_m$ and hence not to the
boundary of $d^*_m$ in view of the relative closeness of $\delta_m$
in $d_m$. Thus, it is sufficient to prove that $\sigma_m\cap\partial
d^*_m\setminus\partial D\neq\varnothing$.

Let us assume that there is a point $x_*\in\sigma_m$ in
$d^*_m\setminus\partial D$. Then there is a point $y_*\in d^*_m$
which is close enough to $\sigma_m$ with $$h(x_0,y_*)\ >\
\frac{1}{2}\left(\rho^-_m+\rho^+_{m+1}\right)$$ because
$h(x_0,y_*)\geqslant\rho^-_m$ and $\rho^+_{m+1}<\rho^-_m$. On the
other hand, there is a point $z_*\in d_{m+1}$ which is close enough
to $\sigma_{m+1}$ such that $$h(x_0,z_*)\ <\
\frac{1}{2}\left(\rho^-_m+\rho^+_{m+1}\right)\,.$$ However, the
points $z_*$ and $y_*$ can be joined by a continuous curve
$\gamma:[0,1]\to d^*_{m+1}$. Note that the sets
$\gamma^{-1}(d^*_m\setminus\overline{d_{m+1}})$ consists of a
countable collection of open disjoint intervals of $[0,1]$ and the
interval $(t_0,1]$ with $t_0\in(0,1)$ and
$z_0=\gamma(t_0)\in\sigma_{m+1}$. Thus, $$h(x_0,z_0)\ <\
\frac{1}{2}\left(\rho^-_m+\rho^+_{m+1}\right)$$ because
$h(x_0,z_0)\leqslant\rho^+_{m+1}$ and $\rho^+_{m+1}<\rho^-_m$. Now,
by the continuity of the function $\varphi(t)=h(x_0,\gamma(t))$,
there is $\tau_0\in(t_0,1)$ such that $$h(x_0,y_0)\ =\
\frac{1}{2}\left(\rho^-_m+\rho^+_{m+1}\right)$$ where
$y_0=\gamma(\tau_0)\in d^*_m$ by the choice of $\gamma$. The
contradiction disproves the above assumption and, thus, the proof is
complete. $\Box$

\bigskip

Later on, given a domain $D$ in ${\mathbb R}^n$, $n\geqslant2$, we
say that a {\bf sequence of points} $x_k\in D$, $k=1,2,\ldots $,
{\bf converges to its end} $K$ if, for every chain $\{ \sigma_m\}$
in $K$ and every domain $d_m$, all points $x_k$ except its finite
collection belong to $d_m$.

\cc
\section{On regular domains}

Recall first of all the following topological notion. A domain
$D\subset{\Bbb R}^n$, $n\geqslant2$, is said to be {\bf locally
connected at a point} $x_0\in\partial D$ if, for every neighborhood
$U$ of the point $x_0$, there is a neighborhood $V\subseteq U$ of
$x_0$ such that $V\cap D$ is connected. Note that every Jordan
domain $D$ in ${\Bbb R}^n$ is locally connected at each point of
$\partial D$, see, e.g., \cite{Wi}, p. 66.

Following \cite{KR} and \cite{KR$_1$}, see also \cite{MRSY} and
\cite{RSal}, we say that $\partial D$ is {\bf weakly flat at a
point} $x_0\in\partial D$ if, for every neighborhood $U$ of the
point $x_0$ and every number $P>0$, there is a neighborhood
$V\subset U$ of $x_0$ such that
\begin{equation}\label{eq1.5KR} M(\Delta(E,F,D))\geqslant P\end{equation} for all continua $E$ and $F$ in $D$
intersecting $\partial U$ and $\partial V$. Here and later on,
$\Delta(E,F,D)$ denotes the family of all paths
$\gamma:[a,b]\to\overline{{\Bbb R}^n}$ connecting $E$ and $F$ in
$D$, i.e., $\gamma(a)\in E$, $\gamma(b)\in F$ and $\gamma(t)\in D$
for all $t\in(a,b)$. We say that the boundary $\partial D$ is {\bf
weakly flat} if it is weakly flat at every point in $\partial D$.

We also say that a {\bf point} $x_0\in\partial D$ is {\bf strongly
accessible} if, for every neighborhood $U$ of the point $x_0$, there
exist a compactum $E$ in $D$, a neighborhood $V\subset U$ of $x_0$
and a number $\delta>0$ such that
\begin{equation}\label{eq1.6KR}M(\Delta(E,F,D))\geqslant\delta\end{equation} for all
continua $F$ in $D$ intersecting $\partial U$ and $\partial V$. We
say that the {\bf boundary} $\partial D$ is {\bf strongly
accessible} if every point $x_0\in\partial D$ is strongly
accessible.

\medskip

\begin{remark}\label{r:13.3}
Here, in the definitions of strongly accessible and weakly flat
boun\-daries, we may take as neighborhoods $U$ and $V$ of a point
$x_0$ only balls (closed or open) centered at $x_0$ or only
neighborhoods of $x_0$ in another fundamental system of
neighborhoods of $x_0$. These concepti\-ons can also be extended in
a natural way to the case of $\overline{{\Bbb R}^n}$ and
$x_0=\infty$. Then we must use the corresponding neighborhoods of
$\infty$.

It is easy to see that if a domain $D$ in ${\Bbb R}^n$ is weakly
flat at a point $x_0\in\partial D$, then the point $x_0$ is strongly
accessible from $D$. Moreover, it was proved by us that if a domain
$D$ in ${\Bbb R}^n$ is weakly flat at a point $x_0\in\partial D$,
then $D$ is locally connected at $x_0$, see, e.g., Lemma 5.1 in
\cite{KR$_1$} or Lemma 3.15 in \cite{MRSY}.
\end{remark}

By the classical geometric definition of V\"ais\"al\"a, see, e.g.,
13.1 in \cite{Va}, a ho\-meo\-mor\-phism $f$ between domains $D$ and
$D'$ in $\Rn$, $n\ge2$, is {\bf
$K$-quasiconformal}\index{$K$-quasiconformal mapping}, abbr. {\bf
$K$-qc mapping}, if
$$
M(\Gamma)/K\ \le\ M(f\Gamma)\ \le\ K\ M(\Gamma)
$$
for every path family $\Gamma$ in $D.$ A homeomorphism $f:D\to D'$
is called {\bf quasi\-con\-for\-mal}\index{quasiconformal mapping},
abbr. {\bf qc}\index{qc mapping}, if $f$ is $K$-quasiconformal for
some $K\in[1,\infty)$, i.e., if the distortion of the moduli of path
families under the mapping $f$ is bounded.

\medskip

We say that the boundary of a domain $D$ in ${\Bbb R}^n$ is {\bf
locally quasiconformal} if every point $x_0\in\partial D$ has a
neighborhood $U$ that can be mapped by a quaisconformal mapping
$\varphi$ onto the unit ball ${\Bbb B}^n\subset{\Bbb R}^n$ in such a
way that $\varphi(\partial D\cap U)$ is the intersection of ${\Bbb
B}^n$ with a coordinate hyperplane. Note that a locally
quasiconformal boundary is weakly flat directly by definitions.

\medskip

In the mapping theory and in the theory of differential equations,
it is often applied the so-called Lipschitz domains whose boundaries
are locally quasiconformal.

\medskip

Recall first that a map $\varphi:X\to Y$ between metric spaces $X$
and $Y$ is said to be {\bf Lipschitz} provided ${\rm
dist}(\varphi(x_1),\varphi(x_2))\leqslant M\cdot{\rm dist}(x_1,x_2)$
for some $M<\infty$ and for all $x_1$ and $x_2\in X$. The map
$\varphi$ is called {\bf bi-Lipschitz} if, in addition, $M^*{\rm
dist}(x_1,x_2)\leqslant{\rm dist}(\varphi(x_1),\varphi(x_2))$ for
some $M^*>0$ and for all $x_1$ and $x_2\in X.$ Later on, $X$ and $Y$
are subsets of ${\Bbb R}^n$ with the Euclidean distance.

\medskip

It is said that a domain $D$ in ${\Bbb R}^n$ is {\bf Lipschitz} if
every point $x_0\in\partial D$ has a neighborhood $U$ that can be
mapped by a bi-Lipschitz homeomorphism $\varphi$ onto the unit ball
${\Bbb B}^n\subset{\Bbb R}^n$ in such a way that $\varphi(\partial
D\cap U)$ is the intersection of ${\Bbb B}^n$ with the a coordinate
hyperplane and $f(x_0)=0$, see, e.g., \cite{Oht}. Note that
bi-Lipschitz homeomorphisms are quasiconformal and hence the
Lipschitz domains have locally quasiconformal boundaries.

\bigskip

We call a bounded domain $D$ in ${\Bbb R}^n$  {\bf regular} if $D$
can be mapped by a quasiconformal mapping onto a domain with locally
quasiconformal boundary.

\medskip

It is clear that each regular domain is finitely connected because
under every homeomorphism between domains $D$ and $D^{\prime}$ in
$\overline{{\Bbb R}^n}$, $n\geqslant2$, there is a natural
one-to-one correspondence between components of the boundaries
$\partial D$ and $\partial D'$, see, e.g., Lemma 5.3 in \cite{IR} or
Lemma 6.5 in \cite{MRSY}. Note also that each finitely connected
domain in the plane whose boundary has no one degenerate component
can be mapped by a conformal mapping onto some domain bounded by a
finite collection of mutually disjoint circles and hence it is a
regular domain, see, e.g., Theorem V.6.2 in \cite{Goluzin}.

\bigskip

As it follows from Theorem 5.1 in \cite{Na}, each prime end of a
regular domain in $\overline{{\Bbb R}^n}$, $n\geqslant2$, is
regular. Combining this fact with Lemma \ref{thabc2} above, we
obtain the following statement.

\medskip

\begin{lemma}\label{thabc3} {\it Each prime end $P$ of a regular
domain $D$ in ${\Bbb R}^n$, $n\geqslant2$, contains a chain of
cross--cuts $\sigma_m$ lying on spheres $S_m$ with center at a point
$x_0\in\partial D$ and with euclidean radii $r_m\to0$ as
$m\to\infty$.}
\end{lemma}

\medskip

\begin{remark}\label{METRIC} As it follows from Theorem 4.1 in \cite{Na},
under a quasiconformal mapping $g$ of a domain $D_0$ with a locally
quasiconformal boundary onto a domain  $D$ in ${\Bbb R}^n$,
$n\geqslant2$, there is a natural one-to-one correspondence between
points of $\partial D_0$ and prime ends of the domain $D$ and,
moreover, the cluster sets $C(g, b)$, $b\in\partial D_0$, coincide
with the impression $I(P)$ of the corresponding prime ends $P$ in
$D$.

If $\overline{D}_P$ is the completion of a regular domain $D$ with
its prime ends and $g_0$ is a quasiconformal mapping of a domain
$D_0$ with a locally quasiconformal boundary onto $D$, then it is
natural to determine in $\overline{D}_p$ a metric
$\rho_0(p_1,p_2)=\left|{\widetilde
{g_0}}^{-1}(p_1)-{\widetilde{g_0}}^{-1}(p_2)\right|$ where
${\widetilde {g_0}}$ is the extension of $g_0$ to $\overline {D_0}$
mentioned above.

If $g_*$ is another quasiconformal mapping of a domain $D_*$ with a
locally quasiconformal boundary onto the domain $D$, then the
corresponding metric
$\rho_*(p_1,p_2)=\left|{\widetilde{g_*}}^{-1}(p_1)-{\widetilde{g_*}}^{-1}(p_2)\right|$
generates the same convergence and, con\-se\-quent\-ly, the same
topology in $\overline{D}_P$ as the metric $\rho_0$ because
$g_0\circ g_*^{-1}$ is a quasiconformal mapping between the domains
$D_*$ and $D_0$ that by Theorem 4.1 in \cite{Na} is extended to a
homeomorphism between $\overline{D_*}$ and $\overline{D_0}$. We
call the given topology in the space $\overline{D}_P$ the {\bf
topology of prime ends}.

This topology can be also described in inner terms of the domain $D$
similarly to Section 9.5 in \cite{CL}, however, we prefer the
definition through the metrics because it is more clear, more
convenient and it is important for us just metrizability of
$\overline{D}_P$. Note also that the space $\overline{D}_P$ for
every regular domain $D$ in ${\Bbb R}^n$ with the given topology is
compact because the closure of the domain $D_0$ with locally
quasiconformal boundary is a compact space and by the construction
$\widetilde {g_0}:\overline{D}_P\to{\overline {D_0}}$ is a
homeomorphism.
\end{remark}

\bigskip

Later on, we will mean the continuity of mappings $f:
\overline{D}_P\to\overline{D^{\prime}}_P$ just with respect to this
topology.

\cc
\section{Continuous extension of lower $Q$-homeomorphisms}

\begin{lemma}\label{l:6.3} {\it Let $D$ and $D'$ be regular domains in
${\Bbb R}^n$, $n\geqslant2$, and $f:D\to D'$ be a lower
$Q$-homeomorphism. If
\begin{equation}\label{e:6.4}\int\limits_{0}^{\delta({x_0})}
\frac{dr}{||\,Q||\,_{n-1}(x_0, r)}\ =\ \infty\qquad\qquad \forall\
x_0\in\partial D\end{equation} for some $\delta({x_0}) <d(x_0) =
\sup\limits_{x\in D}\,|\,x-x_0| $ where
$$\label{e:6.6}||\,Q||\,_{n-1}(x_0,r)\ =\ \left(\int\limits_{D\cap
S(x_0,r)}Q^{n-1} \,d{\cal A}\right)^\frac{1}{n-1},$$ then $f$ can be
extended to a continuous mapping of $\overline{D}_P$ onto
$\overline{D^{\prime}}_P$.}
\end{lemma}

{\bf Proof.} In view of Remark \ref{METRIC}, with no loss of
generality we may assume that the domain $D^{\prime}$ has locally
quasiconformal boundary and
$\overline{D^{\prime}}_P=\overline{D^{\prime}}$. Again by Remark
\ref{METRIC}, namely by metrizability of spaces $\overline{D}_P$ and
$\overline{D^{\prime}}_P$, it suffices to prove that, for each prime
end $P$ of the domain $D$, the cluster set
$$L=C(P,f):=\left\{y\in{{\Bbb
R}^n}:y=\lim\limits_{k\to\infty}f(x_k), x_k\to P, x_k\in D\right\}$$
consists of a single point $y_0\in\partial D^{\prime}$.

Note that $L\neq\varnothing$ by compactness of the set
$\overline{D^{\prime}}$, and it is a subset of $\partial
D^{\prime}$, see, e.g., Proposition 2.5 in \cite{RSal} or
Proposition 13.5 in \cite{MRSY}. Let us assume that there exist at
least two points $y_0$ and $z_0\in L$. Set $U=B(y_0,r_0)$ where
$0<r_0<|y_0-z_0|$.

Let $x_0\in I(P)\subseteq\partial D$ and let $\sigma_k$,
$k=1,2,\ldots\,$, be a chain of cross--cuts of $D$, lying on spheres
$S_k=S(x_0,r_k)$ from Lemma \ref{thabc3}, with the associated
domains $D_k$, $k=1,2,\ldots $. Then there exist points $y_k$ and
$z_k$ in the domains $D_{k}'=f(D_{k})$ such that $|y_0-y_k|<r_0$ and
$|y_0-z_k|>r_0$ and, moreover, $y_k\to y_0$ and $z_k\to z_0$ as
$k\to\infty$. Let $C_k$ a continuous curves joining $y_k$ and $z_k$
in $D_{k}'$. Note that by the construction $\partial U\cap
C_k\neq\varnothing$.

By the condition of strong accessibility of the point $y_0$, see
Remark \ref{r:13.3}, there is a continuum $E\subset D'$ and a number
$\delta>0$ such that
$$M(\Delta(E,C_k;D'))\ \geqslant\ \delta$$ for all large enough $k$.

Without loss of generality, we may assume that the latter condition
holds for all $k=1,2,\ldots$. Note that $C=f^{-1}(E)$ is a compact
subset of $D$ and hence $\varepsilon_0={\rm dist}(x_0,C)>0$. Again,
with no loss of generality, we may assume that $r_k<\varepsilon_0$
for all $k=1,2,\ldots$.

Let $\Gamma_{m}$ be a family of all continuous curves in $D\setminus
D_m$ joining the sphere $S_{0}=S(x_0,\varepsilon_0)$ and
$\overline{\sigma_m}$, $m=1,2,\ldots$. Note that by the construction
$C_k\subset D_k^{\prime}\subset D_m'$ for all $m\leqslant k$ and,
thus, by the principle of minorization
$M(f(\Gamma_{m}))\geqslant\delta$ for all $m=1,2,\ldots$.

On the other hand, the quantity $M(f(\Gamma_{m}))$ is equal to the
capacity of the condenser in $D'$ with facings $\overline{D_{m}'}$
and $\overline{f(D\setminus B_0)}$ where $B_0=B(x_0,\varepsilon_0)$,
see, e.g., \cite{Sh}. Thus, by the principle of minorization and
Theorem 3.13 in \cite{Zi}
$$M(f(\Gamma_{m}))\ \leqslant\
\frac{1}{M^{n-1}(f(\Sigma_{m}))}$$ where $\Sigma_{m}$ is the
collection of all intersections of the domain $D$ and the spheres
$S(x_0,\rho)$, $\rho\in(r_m,\varepsilon_0)$, because
$f(\Sigma_{m})\subset\Sigma(f(S_{m}),f(S_{0}))$ where
$\Sigma(f(S_{m}),f(S_{0}))$ consists of all closed subsets of $D'$
separating $f(S_{m})$ and $f(S_{0})$. Finally, by the condition
(\ref{e:6.4}) we obtain that $M(f(\Gamma_{m}))\to0$ as $m\to\infty$.

The obtained contradiction disproves the assumption that the cluster
set $C(P,f)$ consists of more than one point. $\Box$

\cc
\section{Extension of the inverses of lower $Q$-homeomorphisms}

\begin{lemma}\label{l:9.1} {\it Let $D$ and $D'$ be regular domains in
${\Bbb R}^n$, $n\geqslant 2$, $P_1$ and $P_2$ be different prime
ends of the domain $D$, $f$ be a lower $Q$-homeomorphism of the
domain $D$ onto the domain $D'$, and let $\sigma_m$, $m=1,2,\ldots$,
be a chain of cross--cuts of the prime end $P_1$ from Lemma
\ref{thabc3}, lying on spheres $S(z_1,r_m)$, $z_1\in I(P_1)$, with
associated domains $D_m$. Suppose that the function $Q$ is
integrable in the degree $n-1$ over the surfaces
\begin{equation}\label{INTERSECTION}
D(r)\ =\ \left\{x\in D:|\,x-z_1|=r\right\}\ =\ D\cap
S(z_1,r)\end{equation} for a set $E$ of numbers $r\in(0,d)$ of a
positive linear measure where $d=r_{m_0}$ and where $m_0$ is a
minimal number such that the domain $D_{m_0}$ does not contain
sequences of points converging to $P_2$. If $\partial D'$ is weakly
flat, then}
\end{lemma}\begin{equation}\label{e:9.2}C(P_1,f)\cap C(P_2,f)\ =\
\varnothing.\end{equation}

Note that in view of metrizability of the completion
$\overline{D}_P$ of the domain $D$ with prime ends, see Remark
\ref{METRIC}, the number $m_0$ in Lemma \ref{l:9.1} always exists.

\medskip

{\bf Proof.} Let us choose $\varepsilon\in(0,d)$ such that
$E_0:=\{r\in E:r\in(\varepsilon,d)\}$ has a positive linear measure.
Such a choice is possible in view of subadditivity of the linear
measure and the exhaustion $E = \cup E_m$ where $E_m = \{r\in
E:r\in(1/m,d)\}\,,$ $m=1,2,\ldots $. Note that by Proposition
\ref{prOS2.2}
\begin{equation}\label{e:9.3}M(f(\Sigma_{\varepsilon}))\
>\ 0\end{equation} where $\Sigma_{\varepsilon}$ is the family of all surfaces
$D(r)$, $r\in(\varepsilon,d)$, from (\ref{INTERSECTION}).

Let us assume that $C_1\cap C_2\neq\varnothing$ where
$C_i=C(P_i,f)$, $i=1,2$. By the construction there is $m_1>m_0$ such
that $\sigma_{m_1}$ lies on the sphere $S(z_1,r_{m_1})$ with
$r_{m_1}<\varepsilon$. Let $D_0=D_{m_1}$ and $D_*\subseteq
D\setminus D_{m_0}$ be a domain associated with a chain of
cross--cuts of the prime end  $P_2$. Let $y_0\in C_1\cap C_2$.
Choose $r_0>0$ such that $S(y_0,r_0)\cap f(D_0)\neq\varnothing$ and
$S(y_0,r_0)\cap f(D_*)\neq\varnothing$.

Set $\Gamma=\Gamma(\overline{D_0},\overline{D_*};D)$.
Correspondingly (\ref{e:9.3}), by the principle of minorization and
Theorem 3.13 in \cite{Zi},
\begin{equation}\label{e:9.4}M(f(\Gamma))\ \leqslant\
\frac{1}{M^{n-1}(f(\Sigma_{\varepsilon}))}\ <\ \infty\, .
\end{equation}
Let $M_0>M(f(\Gamma))$ be a finite number. By the condition
$\partial D'$ is weakly flat and hence there is $r_*\in(0,r_0)$ such
that
$$M(\Delta(E,F;D'))\ \geqslant\ M_0$$
for all continua $E$ and $F$ in $D'$ intersecting the spheres
$S(y_0,r_0)$ and $S(y_0,r_*)$. However, these spheres can be joined
by continuous curves $c_1$ and $c_2$ in the domains $f(D_0)$ and
$f(D_*)$ and, in particular, for these curves
\begin{equation}\label{e:9.4a}M_0\ \leqslant\
M(\Delta(c_1,c_2;D'))\ \leqslant\ M(f(\Gamma))\,.\end{equation} The
obtained contradiction disproves the assumption that $C_1\cap
C_2\neq\varnothing$. $\Box$

\medskip

\begin{theorem}\label{t:9.5} {\it Let $D$ and $D'$ be regular domains in ${\Bbb R}^n$,
$n\geqslant 2$. If $f$ is a lower $Q$-homeomorphism $D$ onto $D'$
with $Q\in L^{n-1}(D)$, then $f^{-1}$ can be extended to a
continuous mapping of $\overline{D^{\prime}}_P$ onto
$\overline{D}_P$.}
\end{theorem}

\medskip

{\bf Proof.} By Remark \ref{METRIC}, we may assume with no loss of
generality that $D^{\prime}$ is a circular domain,
$\overline{D^{\prime}}_P=\overline{D^{\prime}}$; $C(y_0,
f^{-1})\ne\varnothing $ for every $y_0\in
\partial D^{\prime}$ because $\overline{D}_P$ is metrizable and compact.
Moreover, $C(y_0, f^{-1})\cap D=\varnothing $, see, e.g.,
Proposition 2.5 in \cite{RSal} or Proposition 13.5 in \cite{MRSY}.

Let us assume that there is at least two different prime ends $P_1$
and $P_2$ in $C(y_0, f^{-1})$. Then $y_0\in C(P_1,f)\cap C(P_2,f)$
and, thus, (\ref{e:9.2}) does not hold. Let $z_1\in\partial D$ be a
point corresponding to $P_1$ from Lemma \ref{thabc3}. Note that
\begin{equation}\label{eqKPR6.2ad} E\ =\ \{r\in(0,\delta):\  Q|_{D\cap S(z_1,r)}\in L^{1}(D\cap S(z_1,r))\}
\end{equation} has a positive linear measure for every $\delta>0$ by the Fubini theorem, see, e.g., \cite{Sa}, because $Q\in
L^{1}(D)$.  The obtained contradiction with Lemma \ref{l:9.1} shows
that $C(y_0, f^{-1})$ contains only one prime end of $D$.

Thus, we have the extension $g$ of $f^{-1}$ to
$\overline{D^{\prime}}$ such that $C(\partial D^{\prime},
f^{-1})\subseteq \overline{D}_P\setminus D$. Really $C(\partial
D^{\prime}, f^{-1})=\overline{D}_P\setminus D$. Indeed, if $P_0$ is
a prime end of $D$, then there is a sequence $x_n$ in $D$ being
convergent to $P_0$. We may assume without loss of generality that
$x_n\to x_0\in\partial D$ and $f(x_n)\to y_0\in\partial D^{\prime}$
because $\overline{D}$ and $\overline{D^{\prime}}$ are compact.
Hence $P_0\in C(y_0, f^{-1})$.

Finally, let us show that the extended mapping
$g:\overline{D^{\prime}}\to\overline{D}_P$ is continuous. Indeed,
let $y_n\to y_0$ in $\overline{D^{\prime}}$. If $y_0\in D^{\prime}$,
then the statement is obvious. If $y_0\in\partial D^{\prime}$, then
take $y^*_n\in D^{\prime}$ such that $|y_n-y^*_n|<1/n$ and
$\rho(g(y_n),g(y^*_n))<1/n$ where $\rho$ is one of the metrics in
Remark \ref{METRIC}. Note that by the construction $g(y^*_n)\to
g(y_0)$ because $y^*_n\to y_0$. Consequently, $g(y_n)\to g(y_0)$,
too. $\Box$

\medskip

\begin{theorem}\label{t:9.12} {\it Let $D$ and $D'$ be regular domains in ${\Bbb
R}^n$, $n\geqslant 2$. If $f:D\to D'$ is a lower $Q$-homeomorphism
with condition (\ref{e:6.4}), then $f^{-1}$ can be extended to a
continuous mapping of $\overline{D'}_P$ onto
$\overline{D}_P$.}\end{theorem}

\medskip

{\bf Proof.} Indeed, by Lemma 9.2 in \cite{KR$_1$} or Lemma 9.6 in
\cite{MRSY}, condition (\ref{e:6.4}) implies that
\begin{equation}\label{e:6.4d}\int\limits_{0}^{\delta}
\frac{dr}{||Q||_{}(x_0,r)}\ =\ \infty\qquad\qquad \forall\
x_0\in\partial D \qquad\forall\
\delta\in(0,\varepsilon_0)\end{equation} and, thus, the set
\begin{equation}\label{eqKPR6.2ad} E\ =\ \{r\in(0,\delta):\  Q|_{D\cap S(x_0,r)}\in L^{1}(D\cap S(x_0,r))\}
\end{equation}
has a positive linear measure for all $x_0\in\partial D$ and all
$\delta\in(0,\varepsilon_0)$ . The rest of arguments is perfectly
similar to one in the proof of the previous theorem. $\Box$
\bigskip

\cc
\section{Homeomorphic extension of lower $Q$-homeomorphisms}

Combining Lemma \ref{l:6.3} and Theorem \ref{t:9.12}, we obtain the
next conclusion.

\medskip

\begin{theorem}\label{t:10.1} {\it Let $D$ and $D'$ be regular domains in ${\Bbb R}^n$,
$n\geqslant 2$, and let $f:D\to D'$ be a lower $Q$-homeomorphism
with
\begin{equation}\label{e:10.2}\int\limits_{0}^{\delta(x_0)}
\frac{dr}{||\,Q||\,_{n-1}(x_0,r)}\ =\ \infty\ \ \ \ \ \ \forall\
x_0\in\partial D\end{equation} for some $\delta(x_0)\in(0,d(x_0))$
where $d(x_0)\ =\ \sup\limits_{x\in D}\,|\,x-x_0|$ and
$$||\,Q||\,_{n-1}(x_0,r)\ =\ \left(\int\limits_{D\cap
S(x_0,r)}Q^{n-1}(x)\,d{\cal A}\right)^\frac{1}{n-1}.$$ Then  $f$ can
be extended to a homeomorphism of $\overline{D}_P$ onto
$\overline{D'}_P$.}
\end{theorem}

\medskip

\begin{corollary}\label{thOSKRSS100} {\it In particular, the conclusion of Theorem
\ref{t:10.1} holds if
\begin{equation}\label{eqOSKRSS100d}q_{x_0}(r)=O\left(\left[\log{\frac1r}\right]^{n-1}\right)\ \ \ \ \ \ \forall\
x_0\in\partial D \end{equation} as $r\to0$ where $q_{x_0}(r)$ is the
mean integral value of $Q^{n-1}$ over the sphere $|x-x_0|=r$. }
\end{corollary}

\medskip

Using Lemma 2.2 in \cite{RS}, see also Lemma 7.4 in \cite{MRSY}, by
Theorem \ref{t:10.1} we obtain the following general lemma that, in
turn, makes possible to obtain new criteria in a great number.

\medskip

\begin{lemma}\label{lemOSKRSS1000} {\it Let $D$ and $D'$ be regular domains in ${\Bbb R}^n$,
$n\geqslant 2$, and let $f:D\to D'$ be a lower $Q$-homeomorphism.
Suppose that
\begin{equation}\label{eqOSKRSS1000}
\int\limits_{D(x_0,\varepsilon)}Q^{n-1}(x)\cdot\psi_{x_0,\varepsilon}^n(|x-x_0|)\,dm(x)=
o\left(I_{x_0}^n(\varepsilon)\right)\qquad\forall x_0\in\partial D
\end{equation} as $\varepsilon\to0$ where $D(x_0,\varepsilon)=\{x\in
D:\varepsilon<|x-x_0|<\varepsilon_0\}$ for
$\varepsilon_0=\varepsilon(x_0)>0$ and where
$\psi_{x_0,\varepsilon}(t): (0,\infty)\to [0,\infty]$,
$\varepsilon\in(0,\varepsilon_0)$, is a two-parameter family of
measurable functions such that
$$0<I_{x_0}(\varepsilon)=
\int\limits_{\varepsilon}^{\varepsilon_0}\psi_{x_0,\varepsilon}(t)\,dt<\infty\qquad\qquad\forall
\varepsilon\in(0,\varepsilon_0)\ .$$ Then  $f$ can be extended to a
homeomorphism of $\overline{D}_P$ onto $\overline{D'}_P$. }
\end{lemma}

\medskip

\begin{remark}\label{rmKR2.9} Note that (\ref{eqOSKRSS1000}) holds, in particular, if
\begin{equation}\label{eqOSKRSS100a} \int\limits_{B(x_0,
\varepsilon_0)}Q^{n-1}(x)\cdot\psi^n
(|x-x_0|)\,dm(x)<\infty\qquad\qquad \forall x_0\in\partial
D\end{equation} where $B(x_0,\varepsilon_0)=\{x\in {\Bbb
R}^n:|x-x_0|<\varepsilon_0\}$ for some
$\varepsilon_0=\varepsilon(x_0)>0$ and where $\psi(t): (0,\infty)\to
[0,\infty]$ is a measurable function such that
$I_{x_0}(\varepsilon)\to\infty$ as $\varepsilon\to0$. In other
words, for the extendability of $f$ to a homeomorphism of
$\overline{D}_P$ onto $\overline{D'}_P$, it suffices the integrals
in (\ref{eqOSKRSS100a}) to be convergent for some nonnegative
function $\psi(t)$ that is locally integrable on $(0,\varepsilon_0]$
but it has a non-integrable singularity at zero.
\end{remark}

Let $D$ be a domain in ${\Bbb R}^n$, $n\geqslant1$. Recall that a
real valued function $\varphi\in L^1_{\rm loc}(D)$ is said to be of
{\bf bounded mean oscillation} in $D$, abbr. $\varphi\in{\rm
BMO}(D)$ or simply $\varphi\in{\rm BMO}$, if
\begin{equation}\label{eq1.11} \Vert\varphi\Vert_*=
\sup\limits_{B\subset D}\ \ \ \dashint\limits_B\vert\varphi(z)-
\varphi_B\vert\,dm(z)<\infty\end{equation} where the supremum is
taken over all balls $B$ in $D$ and
\begin{equation}\label{eq1.12} \varphi_B=\dashint\limits_B
\varphi(z)\,dm(z)=\frac{1}{|B|}\int\limits_B\varphi(z)\,
dm(z)\end{equation} is the mean value of the function $\varphi$ over
$B$. Note that $L^{\infty}(D)\subset{\rm BMO}(D)\subset L^p_{\rm
loc}(D)$ for all $1\leqslant p<\infty$, see, e.g., \cite{ReRy}.

\medskip

A function $\varphi$ in BMO is said to have {\bf vanishing mean
oscillation}, abbr. $\varphi\in{\rm VMO}$, if the supremum in
(\ref{eq1.11}) taken over all balls $B$ in $D$ with
$|B|<\varepsilon$ converges to $0$ as $\varepsilon\to0$. VMO has
been introduced by Sarason in \cite{Sarason}. There are a number of
papers devoted to the study of partial differential equations with
coefficients of the class VMO, see, e.g., \cite{CFL}, \cite{ISbord},
\cite{MRV$^*$}, \cite{Pal} and \cite{Ra}.

\medskip

Following \cite{IR}, we say that a function $\varphi:{\Bbb
R}^n\to{\Bbb R}$, $n\geqslant 2$, has {\bf finite mean oscillation}
at a point $x_0$, write $\varphi\in\ $FMO$(x_0)$, if $\varphi\in
L^1_{\rm loc}$ and
\begin{equation}\label{eq49}
\overline{\lim\limits_{\varepsilon\rightarrow 0}}\ \
\dashint_{B(x_0,\varepsilon)}|\varphi(x)-\widetilde{\varphi}_{\varepsilon}|\,dm(x)<\infty
\end{equation} where $\widetilde{\varphi}_{\varepsilon}$ denotes the mean integral value of the function
$\varphi$ over the ball $B(x_0,\varepsilon)$. We also write
$\varphi\in{\rm FMO}(D)$ or simply $\varphi\in{\rm FMO}$ by context
if this property holds at every point $x_0\in D$. Clearly that ${\rm
BMO}\subset{\rm FMO}$. By definition ${\rm FMO}\subset L^1_{\rm
loc}$ but FMO is not a subset of $L^p_{\rm loc}$ for any $p>1$, see
\cite{MRSY}. Thus, the class FMO is essentially more wide than ${\rm
BMO}_{\rm loc}$.

\medskip

Choosing in Lemma \ref{lemOSKRSS1000} $\psi(t):=\frac{1}{t\log 1/t}$
and applying Corollary 2.3 on FMO in \cite{IR}, see also Corollary
6.3 in \cite{MRSY}, we obtain the next result.

\medskip

\begin{theorem}\label{thOSKRSS101} {\it Let $D$ and $D'$ be regular domains in ${\Bbb R}^n$,
$n\geqslant 2$, and let $f:D\to D'$ be a lower $Q$-homeomorphism. If
$Q^{n-1}(x)$ has finite mean oscillation at every point
$x_0\in\partial D$, then $f$ can be extended to a homeomorphism of
$\overline{D}_P$ onto $\overline{D'}_P$.}
\end{theorem}

\medskip

\begin{corollary}\label{corOSKRSS6.6.2} {\it In particular, the conslusion of Theorem \ref{thOSKRSS101} holds if
\begin{equation}\label{eqOSKRSS6.6.3}
\overline{\lim\limits_{\varepsilon\to0}}\
\dashint_{B(x_0,\varepsilon)}Q^{n-1}(x)\ dm(x)\ <\
\infty\qquad\qquad \forall\ x_0\in\partial D\end{equation}}
\end{corollary}

\medskip

Recall that a point $x_0$ is called a {\bf Lebesgue point} of a
function $\varphi:D\to{\Bbb R}$ if $\varphi$ is integrable in a
neighborhood of $x_0$ and \begin{equation}\label{FMO_eq2.7a}
\lim\limits_{\varepsilon\to 0}\ \ \ \mathchoice
{{\setbox0=\hbox{$\displaystyle{\textstyle -}{\int}$}
\vcenter{\hbox{$\textstyle -$}}\kern-.5\wd0}}
{{\setbox0=\hbox{$\textstyle{\scriptstyle -}{\int}$}
\vcenter{\hbox{$\scriptstyle -$}}\kern-.5\wd0}}
{{\setbox0=\hbox{$\scriptstyle{\scriptscriptstyle -}{\int}$}
\vcenter{\hbox{$\scriptscriptstyle -$}}\kern-.5\wd0}}
{{\setbox0=\hbox{$\scriptscriptstyle{\scriptscriptstyle -}{\int}$}
\vcenter{\hbox{$\scriptscriptstyle -$}}\kern-.5\wd0}}
\!\int_{B(x_0,\varepsilon)}|{\varphi}(x)-{\varphi}(x_0)|\,dm(x)=0\,.
\end{equation}

\medskip

\begin{corollary}\label{corOSKRSS6.6.33} {\it The conslusion of Theorem \ref{thOSKRSS101} holds if
every point $x_0\in\partial D$ is a Lebesgue point of the function
$Q:{{\Bbb R}^n}\to(0,\infty)$.}
\end{corollary}

\medskip

The next statement also follows from Lemma \ref{lemOSKRSS1000} under
the choice $\psi(t)=1/t.$

\medskip

\begin{theorem}\label{thOSKRSS102} {\it Let $D$ and $D'$ be regular domains in ${\Bbb R}^n$,
$n\geqslant 2$, and $f:D\to D'$ be a lower $Q$-homeomorphism. If,
for some $\varepsilon_0=\varepsilon(x_0)>0$, as $\varepsilon\to 0$
\begin{equation}\label{eqOSKRSS10.336a}\int\limits_{\varepsilon<|x-x_0|<\varepsilon_0}Q(x)\,\frac{dm(x)}{|x-x_0|^n}
=o\left(\left[\log\frac{1}{\varepsilon}\right]^n\right)\qquad\qquad
\forall\ x_0\in\partial D\ ,\end{equation}  then $f$ can be extended
to a homeomorphism of $\overline{D}_P$ onto $\overline{D'}_P$.}
\end{theorem}

\medskip

\begin{remark}\label{rmOSKRSS200} Choosing in Lemma \ref{lemOSKRSS1000}
the function $\psi(t)=1/(t\log 1/t)$ instead of $\psi(t)=1/t$,
(\ref{eqOSKRSS10.336a}) can be replaced by the more weak condition
\begin{equation}\label{eqOSKRSS10.336b}\int\limits_{\varepsilon<|x-x_0|<\varepsilon_0}\frac{Q(x)\,dm(x)}{|x-x_0|\,\log{\frac{1}{|x-x_0|}}}
=o\left(\left[\log\log\frac{1}{\varepsilon}\right]^n\right)\end{equation}
and (\ref{eqOSKRSS100d}) by the condition
\begin{equation}\label{eqOSKRSS10.336h} q_{x_0}(r)=o
\left(\left[\log\frac{1}{r}\log\,\log\frac{1}{r}\right]^{n-1}
\right).\end{equation} Of course, we could to give here the whole
scale of the corresponding condition of the logarithmic type using
suitable functions $\psi(t).$
\end{remark}

\medskip

Theorem \ref{t:10.1} has a magnitude of other fine consequences, for
instance:

\medskip

\begin{theorem}\label{thOSKRSS103} {\it Let $D$ and $D'$ be regular domains in ${\Bbb R}^n$,
$n\geqslant 2$, and let $f:D\to D'$ be a lower $Q$-homeomorphism
with
\begin{equation}\label{eqOSKRSS10.36b} \int\limits_D\Phi\left(Q^{n-1}(x)\right)dm(x)<\infty\end{equation}
for a nondecreasing convex function $\Phi:[0,\infty]\to[0,\infty]$
such that, for some $\delta>\Phi(0)$,
\begin{equation}\label{eqOSKRSS10.37b}
\int\limits_{\delta}^{\infty}\frac{d\tau}{\tau\left[\Phi^{-1}(\tau)\right]^{\frac{1}{n-1}}}=
\infty\ .\end{equation} Then $f$ can be extended to a homeomorphism
of $\overline{D}_P$ onto $\overline{D'}_P$.}
\end{theorem}

\medskip

Indeed, by Theorem 3.1 and Corollary 3.2 in
 \cite{RSY}, (\ref{eqOSKRSS10.36b}) and
(\ref{eqOSKRSS10.37b}) imply (\ref{e:10.2}) and, thus, Theorem
\ref{thOSKRSS103} is a direct consequence of Theorem \ref{t:10.1}.

\medskip

\begin{corollary}\label{corOSKRSS6.6.3} {\it In particular, the conclusion of Theorem
\ref{thOSKRSS101} holds if
\begin{equation}\label{eqOSKRSS6.6.6}
\int\limits_{D}e^{\alpha Q^{n-1}(x)}\ dm(x)\ <\ \infty\end{equation}
for some $\alpha>0$.}
\end{corollary}

\medskip

\begin{remark}\label{rmOSKRSS200000}
Note that the condition (\ref{eqOSKRSS10.37b}) is not only
sufficient but also necessary for a cotinuous extension to the
boundary of the mappings $f$ with integral restrictions of the form
(\ref{eqOSKRSS10.36b}), see, e.g., Theorem 5.1 and Remark 5.1 in
\cite{KR$_3$}.

Moreover, by Theorem 2.1 in \cite{RSY}, see also Proposition 2.3 in
\cite{RS1}, (\ref{eqOSKRSS10.37b}) is equivalent to every of the
conditions from the following series:
\begin{equation}\label{eq333Y}\int\limits_{\delta}^{\infty}
H'_{n-1}(t)\ \frac{dt}{t}=\infty\ ,\quad\ \delta>0\ ,\end{equation}
\begin{equation}\label{eq333F}\int\limits_{\delta}^{\infty}
\frac{dH_{n-1}(t)}{t}=\infty\ ,\quad\ \delta>0\ ,\end{equation}
\begin{equation}\label{eq333B}
\int\limits_{\delta}^{\infty}H_{n-1}(t)\ \frac{dt}{t^2}=\infty\
,\quad\ \delta>0\ ,
\end{equation}
\begin{equation}\label{eq333C}
\int\limits_{0}^{\Delta}H_{n-1}\left(\frac{1}{t}\right)\,dt=\infty\
,\quad\ \Delta>0\ ,
\end{equation}
\begin{equation}\label{eq333D}
\int\limits_{\delta_*}^{\infty}\frac{d\eta}{H_{n-1}^{-1}(\eta)}=\infty\
,\quad\ \delta_*>H_{n-1}(+0)\ ,
\end{equation}
\begin{equation}\label{eq333A}
\int\limits_{\delta_*}^{\infty}\,\frac{d\tau}{\tau\Phi_{n-1}^{-1}(\tau)}=\infty\
,\quad\ \delta_*>\Phi(+0)\ ,
\end{equation}
where
\begin{equation}\label{eq333E}
H_{n-1}(t)=\log\Phi_{n-1}(t)\,,\quad\Phi_{n-1}(t)=\Phi\left(t^{n-1}\right)\
.\end{equation}

Here, in (\ref{eq333Y}) and (\ref{eq333F}), we complete the
definition of integrals by $\infty$ if $\Phi_{n-1}(t)=\infty$,
correspondingly, $H_{n-1}(t)=\infty$, for all $t\geqslant T\in{{\Bbb
R}^+}$. The integral in (\ref{eq333F}) is understood as the
Lebesgue--Stieltjes integral and the integrals in (\ref{eq333Y}) and
(\ref{eq333B})--(\ref{eq333A}) as the ordinary Lebesgue integrals.

It is necessary to give one more explanation. From the right hand
sides in the conditions (\ref{eq333Y})--(\ref{eq333A}) we have in
mind $+\infty$. If $\Phi_{n-1}(t)=0$ for $t\in[0,t_*]$, then
$H_{n-1}(t)=-\infty$ for $t\in[0,t_*]$ and we complete the
definition $H_{n-1}'(t)=0$ for $t\in[0,t_*]$. Note, the conditions
(\ref{eq333F}) and (\ref{eq333B}) exclude that $t_*$ belongs to the
interval of integrability because in the contrary case the left hand
sides in (\ref{eq333F}) and (\ref{eq333B}) are either equal to
$-\infty$ or indeterminate. Hence we may assume in
(\ref{eq333Y})--(\ref{eq333C}) that $\delta>t_0$, correspondingly,
$\Delta<1/t_0$ where $t_0:=\sup\limits_{\Phi_{n-1}(t)=0}t$, $t_0=0$
if $\Phi_{n-1}(0)>0$.

\medskip

The most interesting of the above conditions is (\ref{eq333B}) that
can be rewritten in the following form:
\begin{equation}\label{eq5!}
\int\limits_{\delta}^{\infty}\log \Phi(t)\ \
\frac{dt}{t^{n^{\,\prime}}}\ =\ \infty
\end{equation}
where $\frac{1}{n^{\,\prime}}+\frac{1}{n}=1,$ i.e. $n^{\,\prime}=2$
for $n=2,$ $n^{\,\prime}$ is strictly decreasing in $n$ and
$n^{\prime}=n/(n-1)\rightarrow 1$ as $n\rightarrow \infty.$
\end{remark}

\medskip

The theory of the boundary behavior for the lower $Q$-homeomorphisms
de\-ve\-lo\-ped here will find its applications, in particular, to
mappings in classes of Sobolev and Orlicz-Sobolev and also to
finitely bilipschitz mappings that a far--reaching extension of the
well-known classes of isometric and quasiisometric mappings, see,
e.g., \cite{KPR}, \cite{KPRS}, \cite{KRSS1}, \cite{KRSS},
\cite{KSS}, \cite{MRSY} and \cite{RSSY}.

\section{Lower $Q$-homeomorphisms and Orlicz--Sobolev classes}

Following Orlicz, see \cite{Or1}, see also the monographs \cite{KR}
and \cite{Za}, given a convex increasing function $\varphi:{\Bbb
R}^+\to{\Bbb R}^+$, $\varphi(0)=0$, denote by $L^{\varphi}$ the
space of all functions $f:D\to{\Bbb R}$ such that
\begin{equation}\label{eqOS1.1}
\int\limits_{D}\varphi\left(\frac{|f(x)|}{\lambda}\right)\,dm(x)<\infty\end{equation}
for some $\lambda>0$ where $dm(x)$ corresponds to the Lebesgue
measure in $D$. $L^{\varphi}$ is called the {\bf Orlicz space}. In
other words, $L^{\varphi}$ is the cone over the class of all
functions $g:D\to{\Bbb R}$ such that
\begin{equation}\label{eqOS1.2}
\int\limits_{D}\varphi\left(|g(x)|\right)\,dm(x)<\infty\end{equation}
which is also called the {\bf Orlicz class}, see \cite{BO}.

\medskip

The {\bf Orlicz--Sobolev class} $W^{1,\varphi}(D)$ is the class of
all functions $f\in L^1(D)$ with the first distributional
derivatives whose gradient $\nabla f$ belongs to the Orlicz class in
$D$. $f\in W^{1,\varphi}_{\rm loc}(D)$ if $f\in W^{1,\varphi}(D_*)$
for every domain $D_*$ with a compact closure in $D$. Note that by
definition $W^{1,\varphi}_{\rm loc}\subseteq W^{1,1}_{\rm loc}$. As
usual, we write $f\in W^{1,p}_{\rm loc}$ if $\varphi(t)=t^p$,
$p\geqslant1$. Later on, we also write $f\in W^{1,\varphi}_{\rm
loc}$ for a locally integrable vector-function $f=(f_1,\ldots,f_m)$
of $n$ real variables $x_1,\ldots,x_n$ if $f_i\in W^{1,1}_{\rm loc}$
and
\begin{equation}\label{eqOS1.2a} \int\limits_{D}\varphi\left(|\nabla
f(x)|\right)\,dm(x)<\infty\end{equation} where $|\nabla
f(x)|=\sqrt{\sum\limits_{i,j}\left(\frac{\partial f_i}{\partial
x_j}\right)^2}$. Note that in this paper we use the notation
$W^{1,\varphi}_{\rm loc}$ for more general functions $\varphi$ than
in those classic Orlicz classes often giving up the conditions on
convexity and normalization of $\varphi$. Note also that the
Orlicz--Sobolev classes are intensively studied in various aspects
at the moment, see, e.g., \cite{KRSS} and further references
therein.

\medskip

In this connection, recall the minimal definitions which are
relative to Sobolev's classes. Given an open set $U$ in
$\mathbb{R}^n$, $n\ge 2,$ $C_0^{ \infty }(U)$ denotes the collection
of all functions $\psi : U \to \mathbb{R}$ with compact support
having continuous partial derivatives of any order. Now, let $u$ and
$v: U \to {\Bbb R}$ be locally integrable functions. The function
$v$ is called the {\bf distributional derivative} $u_{x_i}$ of $u$
in the variable $x_i$, $i=1,2,\ldots , n$, $x=(x_1,x_2,\ldots ,
x_n)$, if
\begin{equation}\label{eqSTR2.21}
\int\limits_{U} u\, \psi_{x_i} \,dm(x)=- \int\limits_{U} v\, \psi \
dm(x) \ \quad \forall \ \psi \in C_{0}^{\infty}(U)\ .
\end{equation}

\medskip

The {\bf Sobolev classes} $W^{1,p}(U)$ consist of all functions $u:
U \to {\Bbb R}$ in $L^p(U)$ with all distributional derivatives of
the first order in $L^p(U)$. A function $u: U \to {\Bbb R}$ belongs
to $W^{1,p}_{\mathrm{loc}}(U)$ if $u \in W^{1,p}(U_*)$ for every
open set $U_*$ with a compact closure in $ U.$  We use the
abbreviation $W^{1,p}_{\mathrm{loc}}$ if $U$ is either defined by
the context or not essential. The similar notion is introduced for
vector-functions $f: U \to \mathbb{R}^m$ in the component-wise
sense. It is known that a continuous function $f$ belongs to
$W^{1,p}_{\rm loc}$ if and only if $f\in ACL^{p}$, i.e., if $f$ is
locally absolutely continuous on a.e. straight line which is
parallel to a coordinate axis and if the first partial derivatives
of $f$ are locally integrable with the power $p$, see, e.g., 1.1.3
in \cite{Maz}. Recall that the concept of the distributional
(generalized) derivative was introduced by Sobolev in ${\Bbb R}^n$,
$n\geqslant2$, see \cite{So}, and at present it is developed under
wider settings by many authors, see, e.g., many relevant references
in \cite{KRSS}.

\medskip

In this section we show that each homeomorphism $f$ with finite
distortion in ${\Bbb R}^n$, $n\geqslant3$, of the Orlicz--Sobolev
class $W^{1,\varphi}_{\rm loc}$ with the Calderon type condition
\begin{equation}\label{eqOS3.3} \int\limits_{t_*}^{\infty}\left[\frac{t}{\varphi(t)}\right]^
{\frac{1}{n-2}}dt<\infty\end{equation} for some $t_*\in{\Bbb R}^+$,
cf. \cite{Ca}, is a lower $Q$-homeomorphism where $Q=K_f$ is equal
to one of the dilatations of $f$.

Given a mapping $f:D \to \Rn$ with partial derivatives a.e., recall
that $f^\prime(x)$ denotes the Jacobian matrix of $f$ at $x \in D$
if it exists, $J(x)=J(x,f)=\det f^\prime(x)$ is the Jacobian of $f$
at $x$, and $\Vert f^\prime(x)\Vert$ is the operator norm of
$f^\prime(x)$, i.e.,
\begin{equation} \label{eq4.1.2}
\Vert f^\prime(x)\Vert =\max \{|f^\prime(x)h|: h \in \Rn, |h|=1\}.
\end{equation}
We also let
\begin{equation} \label{eq4.1.3}
l(f^\prime(x))= \min \{ |f^\prime(x)h|: h \in \Rn, |h|=1\}.
\end{equation}
The {\bf outer dilatation} of $f$ at $x$ is defined by
\begin{equation} \label{eq4.1.4} K_O(x)=K_O(x,f)= \left
\{\begin{array}{rl}
\frac{\Vert f^\prime(x)\Vert^n}{|J(x,f)|} & {\rm if } \ J(x,f) \neq 0, \\
1 & {\rm if} \ f^\prime(x)=0, \\ \infty & {\rm } \text{otherwise},
\end{array} \right. \end{equation} the {\bf inner dilatation} of $f$ at $x$
by \begin{equation} \label{eq4.1.5} K_I(x)=K_I(x,f)= \left
\{\begin{array}{rl}
\frac{|J(x,f)|}{l(f^\prime(x))^n} &{\rm if } \ J(x,f) \neq 0, \\
1 & {\rm if } \ f^\prime(x)=0, \\
\infty &{\rm} \text{otherwise},
\end{array} \right. \end{equation} Further we also use dilatations $P_O$ and $P_I$  defined by
\begin{equation} \label{eq4.1.44} P_{O}\left( x, f\right)
\,=\,{K^{\frac{1}{n-1}}_{O}(x,f)} \ \ \ \ \  \mbox{and}\ \ \ \ \
P_{I}\left( x, f\right) \,=\,{K^{\frac{1}{n-1}}_{I}(x,f)}\
.\end{equation} Note that
\begin{equation} \label{eq4.1.444} P_O(x,f)\leq K_I(x,f)\ \ \ \ \  \mbox{and}\ \ \ \ \ P_I(x,f) \leq K_O(x,f)\ ,\end{equation} see, e.g.,
Section 1.2.1 in \cite{Re}, and, in particular, $K_O(x,f)$ and
$K_I(x,f),$ $P_O(x,f)$ and $P_I(x,f)$ are simultaneously finite or
infinite. $K_O(x,f) < \infty $ a.e. is equivalent to the condition
that a.e. either $ \det f^\prime(x) > 0$ or $f^\prime(x)=0.$
\medskip

Recall also that a (continuous) mapping $f:D \to \Rn$ is {\bf
absolutely continuous on lines}\index{absolute continuity on lines,
ACL}, abbr. $f\in$ {\bf ACL}, if, for every closed parallelepiped
$P$ in $D$ whose sides are perpendicular to the coordinate axes,
each coordinate function of $f|P$ is absolutely continuous on almost
every line segment in $P$ that is parallel to the coordinate axes.
Note that, if $f\in$ ACL, then $f$ has the first partial derivatives
a.e.
\medskip

In particular, $f$ is ACL if $f\in W^{1,1}_{\mathrm{loc}}$. In
general, mappings in the Sobolev classes
$\bf{W^{1,p}_{\mathrm{loc}}},$\index{$W^{1,p}_{\mathrm{loc}}$} $p\in
[1,\infty ),$ with generalized first partial derivatives in
$L^p_{\mathrm{loc}}$ can be cha\-rac\-te\-ri\-zed as mappings in
$\bf{ACL^p_{\mathrm{loc}}},$ i.e. mappings in ACL whose usual first
partial derivatives are locally integrable in the degree $p;$ see,
e.g., \cite{Maz}, p. 8.
\medskip

\medskip

Now, recall that a homeomorphism $f$ between domains $D$ and $D'$ in
${\Bbb R}^n$, $n\geqslant2$, is called of {\bf finite distortion} if
$f\in W^{1,1}_{\rm loc}$ and \begin{equation}\label{eqOS1.3} \Vert
f'(x)\Vert^n\leqslant K(x)\cdot J_f(x)\end{equation} with some a.e.
finite function $K$. In other words, (\ref{eqOS1.3}) means that
dilatations $K_O(x,f)$, $K_I(x,f),$ $P_O(x,f)$ and $P_I(x,f)$ are
 finite a.e.

\medskip

First this notion was introduced on the plane for $f\in W^{1,2}_{\rm
loc}$ in the work \cite{IS}. Later on, this condition was replaced
by $f\in W^{1,1}_{\rm loc}$ but with the additional condition
$J_f\in L^1_{\rm loc}$ in the monograph \cite{IM}. The theory of the
mappings with finite distortion had many successors, see many
relevant references in the monographs \cite{GRSY} and \cite{MRSY}.
They had as predecessors of the mappings with bounded distortion,
see \cite{Re}, and also\cite{Vo}, in other words, the quasiregular
mappings, see, e.g., \cite{HKM}, \cite{MRV} and \cite{Ri}. They are
also closely connected to the earlier mappings with the bounded
Dirichlet integral and the mappings quasiconformal in the mean which
had a rich history, see, e.g., further references in \cite{MRSY}.

\medskip

Note that the above additional condition $J_f\in L^1_{\rm loc}$ in
the definition of the mappings with finite distortion can be omitted
for homeomorphisms. Indeed, for each homeomorphism $f$ between
domains $D$ and $D'$ in ${\Bbb R}^n$ with the first partial
derivatives a.e. in $D$, there is a set $E$ of the Lebesgue measure
zero such that $f$ satisfies $(N)$-property by Lusin on $D\setminus
E$ and \begin{equation}\label{eqOS1.1.1}
\int\limits_{A}J_f(x)\,dm(x)=|f(A)|\end{equation} for every Borel
set $A\subset D\setminus E$, see, e.g., 3.1.4, 3.1.8 and 3.2.5 in
\cite{Fe}. On this basis, it is also easy by the H\"older inequality
to verify, in particular, that if $f\in W^{1,1}_{\rm loc}$ is a
homeomorphism and $K_f\in L^q_{\rm loc}$ for some $q>n-1$, then also
$f\in W^{1,p}_{\rm loc}$ for some $p>n-1$, that we often use further
to obtain corollaries.

\medskip

On the basis of (\ref{eqOS1.1.1}) below, it is easy to prove the
following useful statement.

\medskip

\begin{proposition}\label{prOS2.4-K} {\it
Let $f$ be an ACL homeomorphism of a domain $D$ in ${\Bbb R}^n$,
$n\geqslant2$, into ${\Bbb R}^n$. Then
$$
(i) \quad  f\in W^{1,1}_{\mathrm loc} \quad  \mbox{if} \quad P_O\in
L^1_{\mathrm loc}\ ,
$$
$$
(ii) \quad  f\in W^{1,\frac{n}{2}}_{\mathrm loc} \quad  \mbox{if}
\quad K_O\in L^1_{\mathrm loc}\ ,
$$
$$
(iii) \quad  f\in W^{1, n-1}_{\mathrm loc} \quad  \mbox{if} \quad
K_O\in L^{n-1}_{\mathrm loc}\ ,
$$
$$
(iv) \quad  f\in W^{1, p}_{\mathrm loc},\ p>n-1 \quad  \mbox{if}
\quad K_O\in L^{{\gamma}}_{\mathrm loc},\ {\gamma}>n-1\ ,
$$
$$
(v) \quad  f\in W^{1, p}_{\mathrm loc},\ p=n{\gamma}/(1+{\gamma})\ge
1 \quad \mbox{if} \quad K_O\in L^{{\gamma}}_{\mathrm loc},\
{\gamma}\ge 1/(n-1)\ .
$$
These conclusions and the estimates (\ref{eqOSadm.1}) are also valid
for all ACL mappings $f: D\to{\Bbb R}^n$ with $J_f\in L^1_{\mathrm
loc}$.}
\end{proposition}

\medskip

Indeed, by the H\"{o}lder inequality applied on a compact set $C$ in
$D$, we obtain on the basis of (\ref{eqOS1.1.1}) the following
estimates of the first partial derivatives
\begin{equation}\label{eqOSadm.1}
\Vert\partial_i f\Vert_p\leqslant\Vert
f^{\prime}\Vert_p\leqslant\Vert K_O^{1/n}\Vert_s\cdot\Vert
J_f^{1/n}\Vert_n\leqslant\Vert
K_O\Vert^{1/n}_{\gamma}\cdot|f(C)|^{1/n}<\infty \end{equation} if
$K_O\in L^{{\gamma}}_{\mathrm loc}$ for some $\gamma\in(0,\infty )$
because $\Vert f^{\prime}(x)\Vert = K_O^{1/n}(x)\cdot J_f^{1/n}(x)$
a.e. where $\frac{1}{p}=\frac{1}{s}+\frac{1}{n}$ and $s={\gamma}n$,
i.e., $\frac{1}{p}=\frac{1}{n}\left(\frac{1}{{\gamma}}+1\right)$.

\medskip

We sometimes use the estimate (\ref{eqOSadm.1}) with no comments to
obtain corollaries.

\medskip

The next statement is key for deriving many consequences of our
theory developed in Sections 5, 6 and 7, cf. Theorem 4.1 in
\cite{KRSS1} and Theorem 5 in \cite{KRSS}.

\medskip

\begin{lemma}\label{thOS4.1} {\it
Let $D$ and $D'$ be domains in ${\Bbb R}^n$, $n\geqslant3$, and let
$\varphi:{\Bbb R}^+\to{\Bbb R}^+$ be a nondecreasing function such
that, for some $t_*\in{\Bbb R}^+$,
\begin{equation}\label{eqOS4.1}
\int\limits_{t_*}^{\infty}\left[\frac{t}{\varphi(t)}\right]^
{\frac{1}{n-2}}dt<\infty\,.\end{equation} Then each homeomorphism
$f:D\to D'$ of finite distortion in the class $W^{1,\varphi}_{\rm
loc}$ is a lower $Q$-homeomorphism at every point
$x_0\in\overline{D}$ with $Q(x)=P_I(x,f)$.} \end{lemma}

\medskip

{\bf Proof.} Let $B$ be a (Borel) set of all points $x\in D$ where
$f$ has a total differential $f'(x)$ and $J_f(x)\ne0$. Then,
applying Kirszbraun's theorem and uni\-que\-ness of approximate
differential, see, e.g., 2.10.43 and 3.1.2 in \cite{Fe}, we see that
$B$ is the union of a countable collection of Borel sets $B_l$,
$l=1,2,\ldots\,$, such that $f_l=f|_{B_l}$ are bi-Lipschitz
homeomorphisms, see, e.g., 3.2.2 as well as 3.1.4 and 3.1.8 in
\cite{Fe}. With no loss of generality, we may assume that the $B_l$
are mutually disjoint. Denote also by $B_*$ the rest of all points
$x\in D$ where $f$ has the total differential but with $f'(x)=0$.

By the construction the set $B_0:=D\setminus\left(B\bigcup
B_*\right)$ has Lebesgue measure zero, see Theorem 1 in \cite{KRSS}.
Hence ${\mathcal A}_S(B_0)=0$ for a.e. hypersurface $S$ in ${\Bbb
R}^n$  and, in particular, for a.e. sphere $S_r:=S(x_0,r)$ centered
at a prescribed point $x_0\in\overline{D}$, see Theorem 2.11 in
\cite{KR$_7$} or Theorem 9.1 in \cite{MRSY}. Thus, by Corollary 4 in
\cite{KRSS} ${\mathcal A}_{S_r^*}(f(B_0))=0$ as well as ${\mathcal
A}_{S_r^*}(f(B_*))=0$ for a.e. $S_r$ where $S_r^*=f(S_r)$.

Let $\Gamma$ be the family of all intersections of the spheres
$S_r$, $r\in(\varepsilon,\varepsilon_0)$,
$\varepsilon_0<d_0=\sup\limits_{x\in D}\,|x-x_0|$, with the domain
$D$. Given $\varrho_*\in{\rm adm}\,f(\Gamma)$ such that
$\varrho_*\equiv0$ outside of $f(D)$, set $\varrho\equiv0$ outside
of $D$ and on $D\setminus B$ and, moreover,
$$
\varrho(x):=\Lambda(x)\cdot\varrho_*(f(x))\qquad{\rm for}\ x\in B
$$
where
$$
\Lambda(x)\ =\ \left[\ J_f(x)\cdot P_I(x,f)\ \right]^{\frac{1}{n}}\
=\ \left[\ \frac{\mbox{det}\, f^{\prime}(x)}{l( f^{\prime}(x))}\
\right]^{\frac{1}{n-1}}\ =$$
$$=\ \left[\
\lambda_2\cdot\ldots\cdot\lambda_n\ \right]^{\frac{1}{n-1}}\
\geqslant \ \left[\ J_{n-1}(x)\ \right]^{\frac{1}{n-1}}\ ;
$$
here as usual $\lambda_n\geqslant\ldots\geqslant\lambda_1$ are
principal dilatation coefficients of $f^{\prime}(x)$, see, e.g.,
Section I.4.1 in \cite{Re}, and $J_{n-1}(x)$ is the
$(n-1)-$dimensional Jacobian of $f|_{S_r}$ at $x$, see Section 3.2.1
in \cite{Fe}.

Arguing piecewise on $B_l$, $l=1,2,\ldots\,$, and taking into
account Kirszbraun's theorem, by Theorem 3.2.5 on the change of
variables  in \cite{Fe}, we have that
$$\int\limits_{S_r}\varrho^{n-1}\,d{\mathcal A}\geqslant\int\limits_{S_{*}^r}\varrho_{*}^{n-1}\,d{\mathcal A}\geqslant1$$
for a.e. $S_r$ and, thus, $\varrho\in{\rm ext\,adm}\,\Gamma$.

The change of variables on each $B_l$, $l=1,2,\ldots\,$, see again
Theorem 3.2.5 in \cite{Fe}, and countable additivity of integrals
give also the estimate
$$\int\limits_{D}\frac{\varrho^n(x)}{P_I(x)}\,dm(x)\leqslant
\int\limits_{f(D)}\varrho^n_*(x)\,dm(x)$$ and the proof is complete.
$\Box$

\begin{corollary}\label{corOS4.1} {\it
Each homeomorphism $f$ with finite distortion in ${\Bbb R}^n$,
$n\geqslant3$, of the class $W^{1,p}_{\rm loc}$ for $p>n-1$ is a
lower $Q$-homeomorphism at every point $x_0\in\overline{D}$ with
$Q=P_I$.}
\end{corollary}

\medskip

Combining the latter and Proposition \ref{prOS2.4-K}, we come to the
following.

\medskip

\begin{corollary}\label{corOS4.1*} {\it
Each homeomorphism $f$ of the class $W^{1,1}_{\rm loc}$ in ${\Bbb
R}^n$, $n\geqslant3$, with $K_O\in L^{q}_{\rm loc}$ for some $q>n-1$
is a lower $Q$-homeomorphism at every point $x_0\in\overline{D}$
with $Q=P_I$.}
\end{corollary}

\medskip

By Proposition \ref{corOS2.1}, we have also the following statement
from Lemma \ref{thOS4.1}.

\medskip

\begin{proposition}\label{corOS2.1os} {\it Let
$f:D\to{\Bbb R}^n$, $n\geqslant3$, be a homeomorphism with $K_I\in
L^1_{\mathrm loc}$ in $W^{1,\varphi}_{\rm loc}$ where $\varphi:{\Bbb
R}^+\to{\Bbb R}^+$ is a nondecreasing function such that
\begin{equation}\label{eqOS4.1os}
\int\limits_{t_*}^{\infty}\left[\frac{t}{\varphi(t)}\right]^
{\frac{1}{n-2}}dt<\infty\,.\end{equation} Then $f$ is a ring
$Q$-homeomorphism at every point $x_0\in\overline{D}$ with $Q=K_I$.}
\end{proposition}

\medskip

\begin{corollary}\label{corOS4.1*os} {\it
Each homeomorphism $f$ of the class $W^{1,1}_{\rm loc}$ in ${\Bbb
R}^n$, $n\geqslant3$, with $K_I\in L^1_{\mathrm loc}$ and $K_O\in
L^{q}_{\rm loc}$ for some $q>n-1$ is a ring $Q$-homeomorphism at
every point $x_0\in\overline{D}$ with $Q=K_I$.}
\end{corollary}

\medskip

\begin{remark}\label{rem7.4.3os} By Remark \ref{rem7.4.3} the conclusion of Proposition
\ref{corOS2.1os} and Corollary \ref{corOS4.1*os} is valid if $K_I$
is integrable only on almost all spheres of small enough radii
centered at $x_0$ assuming that
 the function $K_I$ is extended by zero outside of $D$.
\end{remark}

\section{Boundary behavior of Orlicz--Sobolev classes}

In this section  we assume that $\varphi:{\Bbb R}^+\to{\Bbb R}^+$ is
a nondecreasing function such that, for some $t_*\in{\Bbb R}^+$,
\begin{equation}\label{eqOSKRSS}\int\limits_{t_*}^{\infty}\left[\frac{t}{\varphi(t)}\right]^
{\frac{1}{n-2}}dt<\infty\,.\end{equation} The continuous extension
to the boundary of the inverse mappings has a simpler criterion than
for the direct mappings. Hence we start from the first. Namely, in
view of Lemma \ref{thOS4.1}, we have the following consequence of
Theorem \ref{t:9.5}.

\medskip

\begin{theorem}\label{thKPR8.2} {\it
Let $D$ and $D'$ be regular domains in ${\Bbb R}^n$, $n\geqslant3$
and let $f$ be a homeo\-morphism of $D$ onto $D'$ in a class
$W^{1,\varphi}_{\rm loc}$ with condition (\ref{eqOSKRSS}) and
$K_I\in L^{1}(D)$. Then $f^{-1}$ can be extended to a continuous
mapping of $\overline{D^{\prime}}_P$ onto $\overline{D}_P$.}
\end{theorem}

\medskip

However, as it follows from the example in Proposition 6.3 in
\cite{MRSY}, any degree of integrability $K_I\in L^q(D)$, $q\in
[1,\infty)$, cannot guarantee  the extension by continuity to the
boundary of the direct mappings.

\medskip

By Lemma \ref{thOS4.1}, we have also the following consequence of
Theorem \ref{t:10.1}.

\medskip

\begin{theorem}\label{thKR9.111} {\it
 Let $D$ and $D'$ be regular domains in
${\Bbb R}^n$, $n\geqslant3$, and let $f:D\to D'$ be a homeomorphism
of finite distortion in $W^{1,\varphi}_{\rm loc}$ with condition
(\ref{eqOSKRSS}) such that \begin{equation}\label{eqKPR2.1}
\int\limits_{0}^{\delta(x_0)}\frac{dr}{||K_I||^{\frac{1}{n-1}}(x_0,r)}
=\infty\qquad\forall\ x_0\in\partial D\end{equation} for some
$\delta(x_0)\in(0,d(x_0))$ where $d(x_0)=\sup\limits_{x\in
D}|x-x_0|$ and
$$||K_I||(x_0,r)=\int\limits_{D\cap S(x_0,r)}K_I^{n-1}(x,f)\ d{\mathcal A}\ .$$
Then  $f$ can be extended to a homeomorphism of $\overline{D}_P$
onto $\overline{D'}_P$.} \end{theorem}

\medskip

In particular, as a consequence of Theorem \ref{thKR9.111}, we
obtain the following genera\-lization of the well-known theorems of
Gehring--Martio and Martio--Vuorinen on a homeomorphic extension to
the boundary of quasiconformal mappings between QED domains, see
\cite{GM} and \cite{MV}.

\medskip

\begin{corollary}\label{thKPR9.2} {\it
Let $D$ and $D'$ be regular domains in ${\Bbb R}^n$, $n\geqslant3$,
and let $f:D\to D'$ be a homeo\-morphism of finite distortion in the
class $W^{1,p}_{\rm loc}$, $p>n-1$, in particular,  a
homeo\-morphism in $W^{1,1}_{\rm loc}$ with $K_O\in L^q_{\rm loc}$,
$q>n-1$. If (\ref{eqKPR2.1}) holds, then $f$ can be extended to a
homeomorphism of $\overline{D}_P$ onto $\overline{D'}_P$.}
\end{corollary}

\medskip

By Lemma  \ref{thOS4.1}, as a consequence of Lemma
\ref{lemOSKRSS1000}, we obtain  the following general lemma.

\medskip

\begin{lemma}\label{lemOSKRSS12.1} {\it
 Let $D$ and $D'$ be regular domains in ${\Bbb R}^n$,
$n\geqslant3$, and let $f:D\to D'$ be a homeomorphism of finite
distortion in $W^{1,\varphi}_{\rm loc}$ with condition
(\ref{eqOSKRSS}) such that
\begin{equation}\label{omal} \int\limits_{D(x_0,\varepsilon,\varepsilon_0)}
K_I(x,f)\cdot\psi^n_{x_0,\varepsilon}(|x-x_0|)\,dm(x)=o(I_{x_0}^n(\varepsilon))\
{\rm as}\ \varepsilon\to0\ \forall\ x_0\in\partial D\end{equation}
where $D(x_0,\varepsilon,\varepsilon_0)=\{x\in
D:\varepsilon<|x-x_0|<\varepsilon_0\}$ for some
$\varepsilon_0\in(0,\delta_0)$, $\delta_0=\delta(x_0)=\sup_{x\in
D}|x-x_0|$, and $\psi_{x_0,\varepsilon}(t)$ is a family of
non-negative measurable (by Lebesgue) functions on $(0,\infty)$ such
that \begin{equation}\label{eq5.3} 0\ <\ I_{x_0}(\varepsilon)\ =\
\int\limits_{\varepsilon}^{\varepsilon_0} \psi_{x_0,\varepsilon}(t)\
dt\ <\ \infty\qquad\forall\ \varepsilon\in(0,\varepsilon_0)\
.\end{equation} Then $f$ can be extended to a homeomorphism of
$\overline{D}_P$ onto $\overline{D'}_P$.}
\end{lemma}

\medskip

Choosing in Lemma \ref{lemOSKRSS12.1} $\psi(t)=1/(t\log{1/t})$ and
applying Corollary 2.3 on FMO in \cite{IR}, see also Corollary 6.3
in \cite{MRSY}, we obtain the following result.

\medskip

\begin{theorem}\label{thKPRS12a*} {\it
 Let $D$ and $D'$ be regular domains in
${\Bbb R}^n$, $n\geqslant3$, and let $f:D\to D'$ be a homeomorphism
in $W^{1,\varphi}_{\rm loc}$ with condition (\ref{eqOSKRSS}) such
that
\begin{equation}\label{eqKPRS12a*}K_I(x,f)\ \leqslant\ Q(x)\quad\quad\quad{\rm
a.e.\ in}\ D\end{equation} for a function $Q:{\Bbb R}^n\to{\Bbb
R}^n$, $Q\in{\rm FMO}(x_0)$ for all $x_0\in\partial D$. Then $f$ can
be extended to a homeomorphism of $\overline{D}_P$ onto
$\overline{D'}_P$.}
\end{theorem}

\medskip

In the next consequences, we assume that $K_I(x,f)$ is extended by
zero outside of $D$.

\medskip

\begin{corollary}\label{corKPRS12a*} {\it
 In particular, the conclusions
of Theorem \ref{thKPRS12a*} hold if
\begin{equation}\label{eqKPRS12b*}\overline{\lim\limits_{\varepsilon\to0}}\
\ \dashint_{B(x_0,\varepsilon)}K_I(x,f)\ dm(x)\ <\
\infty\qquad\forall\ x_0\in\partial D\ .\end{equation}}
\end{corollary}

\medskip

Similarly, choosing in Lemma \ref{lemOSKRSS12.1} the function
$\psi(t)=1/t$, we come to the following more general statement.

\medskip

\begin{theorem}\label{thKPRS12b*} {\it
Let $D$ and $D'$ be regular domains in ${\Bbb R}^n$, $n\geqslant3$,
and let $f:D\to D'$ be a homeomorphism in $W^{1,\varphi}_{\rm loc}$
with condition (\ref{eqOSKRSS}) such that
\begin{equation}\label{eqKPRS12c*}
\int\limits_{\varepsilon<|x-x_0|<\varepsilon_0}K_I(x,f)\
\frac{dm(x)}{|x-x_0|^n}\ =\
o\left(\left[\log\frac{\varepsilon_0}{\varepsilon}\right]^n\right)\qquad\forall\
x_0\in\partial D\end{equation} as $\varepsilon\to0$ for some
$\varepsilon_0\in(0,\delta_0)$ where
$\delta_0=\delta(x_0)=\sup_{x\in D}|x-x_0|$. Then $f$ can be
extended to a homeomorphism of $\overline{D}_P$ onto
$\overline{D'}_P$.}
\end{theorem}

\medskip

\begin{corollary}\label{corKPRS12b*} {\it
 The condition (\ref{eqKPRS12c*}) and the
assertion of Theorem \ref{thKPRS12b*} hold if
\begin{equation}\label{eqKPRS12d*}K_I(x,f)\ =\ o\left(\left[\log\frac{1}{|x-x_0|}\right]^{n-1}\right)\end{equation}
as $x\to x_0$. The same holds if \begin{equation}\label{eqKPRS12e*}
k_f(r)=o\left(\left[\log\frac{1}{r}\right]^{n-1}\right)
\end{equation} as $r\to0$ where $k_f(r)$ is the mean value of the
function $K_I(x,f)$ over the sphere $|x-x_0|=r$.}
\end{corollary}

\medskip

\begin{remark}\label{rmKRRSa*} Choosing in Lemma \ref{lemOSKRSS12.1} the function
$\psi(t)=1/(t\log{1/t})$ instead of $\psi(t)=1/t$, we are able to
replace (\ref{eqKPRS12c*}) by
\begin{equation}\label{eqKPRS12f*}
\int\limits_{\varepsilon<|x-x_0|<1}\frac{K_I(x,f)\
dm(x)}{\left(|x-x_0|\log{\frac{1}{|x-x_0|}}\right)^n}
=o\left(\left[\log\log\frac{1}{\varepsilon}\right]^n\right)\end{equation}
and (\ref{eqKPRS12e*}) by
\begin{equation}\label{eqKPRS12g*}
k_f(r)=o\left(\left[\log\frac{1}{r}\log\log\frac{1}{r}\right]^{n-1}\right).\end{equation}
Thus, it is sufficient to require that
\begin{equation}\label{eqKPRS12h*}k_f(r)=O\left(\left[\log\frac{1}{r}\right]^{n-1}\right)
\end{equation}

In general, we could give here the whole scale of the corresponding
conditions in terms of $\log$ using functions $\psi(t)$ of the form
$1/(t\,\log\ldots\log1/t)$. \end{remark}

\medskip

\begin{theorem}\label{thKR4.1} {\it
Let $D$ and $D'$ be regular domains in ${\Bbb R}^n$, $n\geqslant3$,
and let $f:D\to D'$ be a homeomorphism in $W^{1,\varphi}_{\rm loc}$
with condition (\ref{eqOSKRSS}) such that
\begin{equation}\label{eqKR4.1}
\int\limits_{D}\Phi(K_I(x,f))\ dm(x)\ <\ \infty\end{equation} for a
non-decreasing convex function $\Phi:\overline{{\Bbb
R}^+}\to\overline{{\Bbb R}^+}$. If, for some $\delta>\Phi(0)$,
\begin{equation}\label{eqKR4.2}\int\limits_{\delta}^{\infty}\frac{d\tau}{\tau\left[\Phi^{-1}(\tau)\right]^{\frac{1}{n-1}}}=
\infty\end{equation} then $f$ can be extended to a homeomorphism of
$\overline{D}_P$ onto $\overline{D'}_P$.}
\end{theorem}

\medskip

Indeed, by Theorem 3.1 and Corollary 3.2 in
 \cite{RSY}, (\ref{eqKR4.1}) and
(\ref{eqKR4.2}) imply (\ref{eqKPR2.1}) and, thus, Theorem
\ref{thKR4.1} is a direct consequence of Theorem \ref{thKR9.111}.

\medskip

\begin{corollary}\label{corOSKRSS4.6.3} {\it The conclusion of Theorem
\ref{thKR4.1} holds if, for some $\alpha>0$,
\begin{equation}\label{eqOSKRSS4.6.6}
\int\limits_{D}e^{\alpha K_I(x,f)}\ dm(x)\ <\ \infty\
.\end{equation} }
\end{corollary}

\medskip

\begin{remark}\label{rmKR4.1} Note that by Theorem 5.1 and Remark 5.1 in
\cite{KR$_3$} the conditions (\ref{eqKR4.2}) are not only sufficient
but also necessary for continuous extension to the boundary of $f$
with the integral constraint (\ref{eqKR4.1}).

\medskip

Recall that by Remark \ref{rmOSKRSS200000} the condition
(\ref{eqKR4.2}) is equivalent to each of the conditions
(\ref{eq333Y})--(\ref{eq333A}) and, in particular, to the following
condition
\begin{equation}\label{eqKR4.4}
\int\limits_{\delta}^{\infty}\log\,\Phi(t)\,\frac{dt}{t^{n'}}=+\infty\end{equation}
for some $\delta>0$ where $\frac{1}{n'}+\frac{1}{n}=1$, i.e., $n'=2$
for $n=2$, $n'$ is strictly decreasing in $n$ and $n'=n/(n-1)\to1$
as $n\to\infty$. \end{remark}

\medskip

Finally, note that all these results hold, for instance, if $f\in
W^{1,p}_{\rm loc}$, $p>n-1$, and, in particular, if $f\in
W^{1,1}_{\rm loc}$ and $K_O\in L^{q}_{\rm loc}$, $q>n-1$. Moreover,
the results can be extended to Riemannian manifolds, see, e.g.,
\cite{ARS} and \cite{KSS}.
\medskip

\section{On finitely bi--Lipschitz mappings}

Given an open set $\Omega\subseteq\Rn$, $n\geqslant 2$, following
Section 5 in \cite{KR$_7$}, see also Section 10.6 in \cite{MRSY}, we
say that a mapping $f:\Omega\to\Rn$ is {\bf finitely bi-Lipschitz}
if
\begin{equation}\label{eq8.12.2} 0\ <\ l(x,f)\ \leqslant\ L(x,f)\ <\
\infty \ \ \ \ \ \forall\ x\in\Omega\end{equation} where
\begin{equation} \label{eq8.1.6} L(x,f)\ =\
\limsup_{y\to x}\ \frac{|f(y)-f(x)|}{|y-x|}
\end{equation} and
\begin{equation}\label{eq8.1.7} l(x,f)\ =\ \liminf_{y\to
x}\ \frac{|f(y)-f(x)|}{|y-x|}\ ,\end{equation} cf. Section 4 above
for the definition of bi-Lipschitz mappings.

\medskip

By the classic Stepanov theorem, see \cite{Step}, see also
\cite{Maly}, we obtain from the right hand inequality in
(\ref{eq8.12.2}) that finitely bi-Lipschitz mappings are
differentiable a.e. and from the left hand inequality in
(\ref{eq8.12.2}) that $J_f(x)\ne 0$ a.e. Moreover, such mappings
have $(N)-$property with respect to each Hausdorff measure, see,
e.g., either Lemma 5.3 in \cite{KR$_7$} or Lemma 10.6 \cite{MRSY}.
Thus, the proof of the following lemma is perfectly similar to one
of Lemma \ref{thOS4.1} and hence we omit it, cf. also similar but
weaker Corollary 5.15 in \cite{KR$_7$} and Corollary 10.10 in
\cite{MRSY}.

\medskip

\begin{lemma}\label{pr8.12.15} {\it Every finitely bi-Lipschitz
homeomorphism $f:\Omega\to\Rn$, $n\geqslant 2$, is a lower
$Q$-homeomorphism with $Q=P_I$.} \end{lemma}

\medskip

By Proposition \ref{corOS2.1}, we have also the following statement
from Lemma \ref{pr8.12.15}.

\medskip

\begin{proposition}\label{corOS2.1fb} {\it Every finitely bi-Lipschitz
homeomorphism $f:\Omega\to\Rn$, $n\geqslant 2$, with $K_I\in
L^1_{\mathrm loc}$ is a ring $Q$-homeomorphism at each point $x_0\in
\overline{D}$ with $Q=K_I$.}
\end{proposition}

\medskip

\begin{remark}\label{rem7.4.3fb} By Remark \ref{rem7.4.3} the conclusion of Proposition \ref{corOS2.1fb}
is valid if $K_I$ is integrable only on almost all spheres of small
enough radii centered at $x_0$ assuming that
 the function $K_I$ is extended by zero outside of $D$.
\end{remark}

\medskip

\begin{corollary}\label{cor8.12.15} {\it All results on lower $Q-$homeomorphisms
in Sections 5, 6 and 7 are valid for finitely bi-Lipschitz
homeomorphisms $f:\Omega\to\Rn$, $n\ge 2$, with $Q=P_I$.}
\end{corollary}

\medskip

All these results for finitely bi-Lipschitz homeomorphisms are
perfectly si\-mi\-lar to the corresponding results for
homeomorphisms with finite distortion in the Orlich--Sobolev classes
from Section 9. Hence we will not formulate all them in the explicit
form here in terms of inner dilatation $K_I$.

\medskip

We give here for instance only one of these results.

\medskip

\begin{theorem}\label{thKR4.1fb} {\it
Let $D$ and $D'$ be regular domains in ${\Bbb R}^n$, $n\geqslant2$,
and let $f:D\to D'$ be a finitely bi-Lipschitz homeomorphism such
that
\begin{equation}\label{eqKR4.1fb}
\int\limits_{D}\Phi(K_I(x,f))\ dm(x)\ <\ \infty\end{equation} for a
non-decreasing convex function $\Phi:\overline{{\Bbb
R}^+}\to\overline{{\Bbb R}^+}$. If, for some $\delta>\Phi(0)$,
\begin{equation}\label{eqKR4.2fb}\int\limits_{\delta}^{\infty}\frac{d\tau}{\tau\left[\Phi^{-1}(\tau)\right]^{\frac{1}{n-1}}}=
\infty\end{equation} then $f$ can be extended to a homeomorphism of
$\overline{D}_P$ onto $\overline{D'}_P$.}
\end{theorem}

\medskip

\begin{corollary}\label{corOSKRSS4.6.3fb} {\it The conclusion of Theorem
\ref{thKR4.1fb} holds if, for some $\alpha>0$,
\begin{equation}\label{eqOSKRSS4.6.6fb}
\int\limits_{D}e^{\alpha K_I(x,f)}\ dm(x)\ <\ \infty\
.\end{equation} }
\end{corollary}

\medskip
\noindent
{\bf Denis Kovtonyuk and Vladimir Ryazanov,}\\
Institute of Applied Mathematics and Mechanics,\\
National Academy of Sciences of Ukraine,\\
74 Roze Luxemburg Str., Donetsk, 83114, Ukraine,\\
denis$\underline{\ \ }$\,kovtonyuk@bk.ru, vl.ryazanov1@gmail.com

\end{document}